\documentclass[letterpaper, conference]{ieeetran}
\IEEEoverridecommandlockouts
\usepackage{graphics} 
\usepackage{epsfig} 
\usepackage{amsmath} 
\usepackage{amssymb}  
\usepackage{cite}
\usepackage{color}
\usepackage{graphicx}
\usepackage{algpseudocode}
\usepackage{algorithm}
\usepackage{tikz-cd}
\usepackage{adjustbox}

\let\labelindent\relax
\usepackage{enumitem}
\renewcommand{\theenumi}{(\alph{enumi})}

\usepackage[font=footnotesize]{caption}

\newtheorem{theorem}{Theorem}[section]
\newtheorem{definition}[theorem]{Definition}

\newtheorem{lemma}[theorem]{Lemma}

\newtheorem{remark}[theorem]{Remark}
\newtheorem{example}[theorem]{Example}
\newtheorem{corollary}[theorem]{Corollary}
\newtheorem{proposition}[theorem]{Proposition}  
\newtheorem{problem}{Problem}

\newcommand{\naturals}{\mathbb{N}}
\newcommand{\real}{\mathbb{R}}
\newcommand{\cplx}{\mathbb{C}}

\newcommand{\Ac}{\mathcal{A}}

\newcommand{\Ec}{\mathcal{E}}
\newcommand{\Fc}{\mathcal{F}}

\newcommand{\Kc}{\mathcal{K}}

\newcommand{\Rc}{\mathcal{R}}
\newcommand{\Sc}{\mathcal{S}}
\newcommand{\Tc}{\mathcal{T}}
\newcommand{\Uc}{\mathcal{U}}

\newcommand{\Xc}{\mathcal{X}}

\newcommand{\until}[1]{\{1,\dots,#1\}}

\newcommand{\Span}{\operatorname{span}}

\newcommand{\identity}{\operatorname{id}}
\newcommand{\card}{\operatorname{card}}

\newcommand{\restr}[2]{#1 \!\! \restriction_{#2}}

\newcommand{\aug}{\operatorname{aug}}
\newcommand{\Faug}{\Fc^{\aug}}
\newcommand{\Finf}{\Fc^{\infty}}
\newcommand{\ci}{\operatorname{CI}}
\newcommand{\Faugci}{\Faug_{\ci}}
\newcommand{\Finfci}{\Finf_{\ci}}
\newcommand{\Taug}{\Tc^{\aug}}
\newcommand{\Kaug}{\Kc^{\aug}}
\newcommand{\Kinf}{\Kc^{\infty}}

\newcommand{\longthmtitle}[1]{\mbox{}{\textit{(#1):}}}
\newcommand{\setdef}[2]{\{#1 \; | \; #2\}}

\newcommand{\oprocendsymbol}{\hbox{$\square$}}
\newcommand{\oprocend}{\relax\ifmmode\else\unskip\hfill\fi\oprocendsymbol}
\def\eqoprocend{\tag*{$\square$}}

\newcommand{\Rinftoaug}{R_{\Xc \times \ell(\Uc)}^{\Xc \times \Uc}}
\newcommand{\Raugtoorig}{R_{\Xc \times \Uc}^{\Xc}}
\newcommand{\Eaugtoinf}{E_{\Xc \times \Uc}^{\Xc \times \ell(\Uc)}}
\newcommand{\Eorigtoaug}[1]{E_{\Xc, u \equiv #1}^{\Xc \times \Uc}}
\newcommand{\RopFaugtoF}[1]{\Rc_{\Faug}^{\Fc,u \equiv #1}}
\newcommand{\EopFtoFaug}{\Ec_{\Fc}^{\Faug}}
\newcommand{\RopFinftoFaug}{\Rc_{\Finf}^{\Faug}}
\newcommand{\EopFaugtoFinf}{\Ec_{\Faug}^{\Finf}}
\newcommand{\RopCIFinftoFaug}{\Rc_{\Finfci}^{\Faugci}}
\newcommand{\EopCIFaugtoFinf}{\Ec_{\Faugci}^{\Finfci}}
\newcommand{\RopCIFaugtoF}{\Rc_{\Faugci}^{\Fc}}
\newcommand{\EopCIFtoFaug}{\Ec_{\Fc}^{\Faugci}}

\graphicspath{{figures/}}

\renewcommand{\baselinestretch}{1}

\allowdisplaybreaks

\title{\Large
\bf 
Two Roads to Koopman Operator Theory for Control:
\\
Infinite Input Sequences and Operator Families \thanks{A preliminary
version of this work will appear as~\cite{MH-IM-JC:25-cdc} at the
2025 IEEE Conference on Decision and Control. The work of M. Haseli and J. Cort\'es was supported by ONR Award N00014-23-1-2353. The work of I. Mezi\'c was supported by AFOSR Award FA9550-22-1-0531 to AIMdyn, Inc and DARPA under Agreement No. HR00112590152. Approved
for public release; distribution is unlimited.}}

\author{Masih Haseli \quad Igor Mezi\'c \quad Jorge
Cort\'es\thanks{M. Haseli is with the Department of Computing and Mathematical Sciences, California Institute of Technology, mhaseli@caltech.edu. During the bulk of this work, M. Haseli was with the Department of Mechanical and Aerospace Engineering, University of California, San Diego.  I. Mezi\'c is with the
Department of Mechanical Engineering, University of California,
Santa Barbara, mezic@ucsb.edu. J. Cort\'es is with the Department of
Mechanical and Aerospace Engineering, University of California,
San Diego, cortes@ucsd.edu.}}

\begin{document}

\maketitle

\begin{abstract}
The Koopman operator, originally defined for dynamical systems
without input, has inspired many applications in control. Yet, the
theoretical foundations underpinning this progress in control remain
underdeveloped.  This paper investigates the theoretical structure
and connections between two extensions of Koopman theory to control:
(i) Koopman operator via infinite input sequences and (ii) the
Koopman control family. Although these frameworks encode system
information in fundamentally different ways, we show that under
certain conditions on the function spaces they operate on, they are
equivalent.  The equivalence is both in terms of the actions of the
Koopman-based formulations in each framework
as well as the function values on the system trajectories.  Our
analysis provides constructive tools to translate between the
frameworks, offering a unified perspective for Koopman methods in
control.
\end{abstract}

\section{Introduction}
The Koopman operator provides a representation of the evolution of
nonlinear systems through linear operators acting on vector spaces of
functions. This viewpoint enables the use of both the algebraic
structure of vector spaces and the spectral properties of linear
operators to analyze the behavior of nonlinear dynamics. These
features have motivated extensive research on the theoretical
foundations and practical applications of the Koopman operator. The
original formulation was developed for systems without inputs, where
both theory and applications are now well established. In contrast,
for control systems, although the practical side has received
significant attention, the theoretical foundation remains in its early
stages. In this paper, we investigate the connections between two
established theoretical frameworks that extend Koopman operator theory
to control systems. Our analysis studies how each framework encodes the information about the system and how information from one framework can be translated into the other, and establishes conditions on the respective function spaces under which both frameworks are equivalent.

\subsection*{Literature Review}
The Koopman operator represents the dynamics of a nonlinear system
through a linear operator acting on a vector space of
functions~\cite{BOK:31,BOK-JVN:32}. The linearity of both the function
space and the operator lead to well-structured algebraic properties
that can be leveraged for the systematic study of complex nonlinear
systems, particularly in situations where geometric methods are
difficult to apply~\cite{IM:21}.  These properties have sparked
significant research in analyzing complex systems with a myriad of
applications, including stability
analysis~\cite{AM-IM:16,SAD-DVD:23,CMZ-AM:25,JH-EY:21}, signal
processing~\cite{ZZ-JL-YY:24}, fluid
dynamics~\cite{IM:13,AD-DL-DS-AT:21}, power
networks~\cite{YS-IM-TH-FR:16,MN-LM:18,SPN-SG-SK-SP-KA-YW-SC:20},
biological systems~\cite{AH-NB-EY:19,JH-EY:21}, and hybrid
systems~\cite{NK-YS:24}.

While the Koopman operator was initially formulated for systems
without input, the literature has adapted Koopman-inspired methods to
control applications.  Many Koopman-inspired works in control do not
propose a formal extension of the Koopman theory to systems with
input, and instead rely on the idea of \emph{``lifting''} to higher
dimensions inspired by finite-dimensional linear forms associated with
the Koopman operator for systems without input.  These lifting
techniques, combined with methods from classical system theory and
differential geometry, have proven effective in practice. Among them,
\emph{lifted linear} models~\cite{JLP-SLB-JNK:16} are the simplest and
most widely used, as they allow to leverage highly efficient linear
control techniques such as linear quadratic regulators (LQR) and
linear model predictive control (MPC). While lifted linear models lack
the structural richness needed to capture cross terms between inputs
and states or input nonlinearities, see
e.g.,~\cite{RS-KW-IM-JB-MS-FG:25} for a discussion, they remain
effective in many applications, since feedback loops or MPC schemes
can often compensate for the resulting model mismatch.  The
work~\cite{DU-KD:25} addresses some of the limitations of lifted
linear models via a two-stage learning scheme using orthogonal and
oblique projections and~\cite{HHA-JASC:24} provides a method to find
lifted linear models based on the physical system structure.  For the
special case of control-affine systems, one can formulate an
operator-theoretic framework using a family of Koopman generators that
are affine in the input, see e.g.,~\cite{SP-SEO-CWR:20} for a
discussion. This perspective has motivated the development of
finite-dimensional lifted bilinear models~\cite{DG-DAP:21}.  The
work~\cite{RS-MS-JB-KW-FA:25} establishes error bounds for bilinear
models and~\cite{RS-MS-KW-JB-FA:25} provides controllers with
closed-loop stability guarantees. We refer the reader to the recent
overview~\cite{RS-KW-IM-JB-MS-FG:25} and the references therein for a
comprehensive discussion of bilinear models and their properties.  A
different and effective way of modeling control systems that are not
necessarily control-affine via the Koopman operator is to fix finitely
many constant values of the input and control the system by switching
between the constant-input systems. This idea has been explored in
different forms in the literature. The
works~\cite{AS-DE:17,AS-AM-DE:18} have used variation of pulse-based
control for monotone systems. The work~\cite{SP-SK:19} uses
finite-dimensional approximations of the Koopman operators associated
with systems created by setting the input to be constant and turns the
control problem into a switched linear form. The
work~\cite{JH-KC:24-lcss} takes a similar approach and solves the
resulting optimal control problem for switched linear systems via a
Markov Chain Monte Carlo (MCMC) method.  The aforementioned lifted
linear, bilinear, and switch linear models have found their way into
many applications in controls, including optimal
control~\cite{AS-AM-DE:18,UV-DTC:23,VNF-SAD-DVD:25}, feedback
linearization~\cite{DG-VK-FP:24}, safety and reachability
analysis~\cite{CF-YC-ADM-JWB:20,MB-DP:23},
MPC~\cite{MK-IM-automatica:18,XZ-WP-RS-SY-XX:22}, control design
based on control Lyapunov functions~\cite{VZ-EB:23-neurocomputing} and a wide range of robotics applications~\cite{CF-YC-ADM-JWB:20,DAH-MJB-EK-ABC-PCC-IM-EWH:23,HHA-JASC:24,LS-MH-GM-DB-IA-TM-JC-KK:24}.

Despite the widespread use of finite-dimensional Koopman-inspired
lifted models in control and robotics, the theoretical basis of
Koopman extensions for control has remained relatively unexplored.
The work~\cite{MK-IM-automatica:18} formally extends Koopman operator
theory to control systems by considering the all behaviors generated
under all possible infinite input sequences. Specifically, by
augmenting the state space with the space of infinite input sequences,
one constructs a dynamical system (without input) that encapsulates
all trajectories the original control system can produce. A Koopman
operator is then naturally associated with this extended dynamical
system.  Importantly, this representation does not impose restrictive
assumptions on the function space, nor does it depend on particular
structural properties of the system (such as control-affine dynamics).
The work~\cite{MH-JC:25-auto} provides a different extension of
Koopman theory to general (not necessarily input affine) control
systems, termed Koopman Control Family (KCF). This framework
characterizes the system behavior through the collection of Koopman
operators associated with the systems without input obtained by fixing
the control input to a constant value, for all possible values. KCF
framework also provides a finite-dimensional form termed ``input-state
separable'' model which captures the lifted linear, bilinear, and
switched linear models mentioned above as special cases. KCF also
provides methods to compute the optimal approximation of input-state
separable forms as well as error bounds for the prediction of all
functions in the subspace via the notion of invariance
proximity~\cite{MH-JC:23-csl,MH-JC:24-scl}.  Finally, the
work~\cite{ML:25} provides an alternative operator-theoretic approach
to encode the system behavior. This framework relies on a product
space to distinguish between the effect of states and inputs, thus
providing an effective way to capture and predict the system's
behavior accompanied by data-driven methods.

\subsection*{Statement of Contributions}
We study formal extensions of Koopman operator theory to discrete-time
control systems that are not necessarily in control-affine form.  Our
exposition begins by studying two intuitive approaches to extending
Koopman theory to control systems and showing that they are unable to
capture the system's evolution beyond a single time step. This
observations indicates that, contrary to common belief, extending the
Koopman theory to control systems is far from trivial and requires
rigorous theoretical analysis.  We then turn our attention to two
formal extensions of Koopman theory to control systems: (i) the
Koopman operator via infinite input
sequences~\cite{MK-IM-automatica:18} and (ii) the Koopman control
family (KCF)~\cite{MH-JC:25-auto}. These extensions encode system
information in fundamentally different ways: the former employs a
single operator but requires infinite input sequences, whereas the
latter consists of infinitely many operators defined for constant
inputs. To connect these frameworks, we utilize a parameterization of
KCF using a single operator.  Since different paradigms are based on
substantially different function spaces, we establish connections
between them through linear composition operators. With these
connections at hand, we derive precise algebraic relationships between
the frameworks, enabling the translation of information across
them. Notably, under mild conditions on function spaces, we show that
the frameworks are equivalent and capture the same information about
the systems through linear operators. Moreover, we introduce a
different type of equivalence, concerning the evolution of function
\emph{values} along system trajectories. Our analysis shows that the frameworks, albeit structurally different, are two sides of the same coin in terms of representing the control system. Moreover, the algebraic tools we present allow to translate information from one framework to the other, thus unifying them  and providing a consistent Koopman-based representation for control systems.

\subsection*{Notation}
We use $\naturals$, $\naturals_0$, $\real$, and $\cplx$ to represent
natural, non-negative integer, real, and complex numbers.  Given sets
$A$ and $B$, $A \subseteq (\subset) B$ means that $A$ is a (proper)
subset of $B$. Moreover, $A \cup B$ and $A \cap B$ are the union and
intersection of $A$ and $B$. In addition,
$A \times B := \setdef{(a,b)}{a \in A, \; b \in B}$ denotes the
Cartesian product of $A$ and $B$.  We also denote the cardinality of
$A$ by $\card(A)$.  Given a function $f:A \to B$ and the set
$S \subseteq A$, we denote by $\restr{f}{S}: S \to B$, the function
created by restricting the domain\footnote{We also use a similar
notation for restricting the domain of functions from product spaces
$A \times B$ to one of the sets in the pair (e.g., $A$ or $B$). In
such cases the notation is explained where it is used.} of $f$ to
$S$. The image of set $E \subseteq A$ under function $f:A \to B$ is
$f(E):= \setdef{f(y)}{y \in E}$.  Given that we study functions with
various domains and co-domains throughout the paper, when convenient
we use the domain and co-domain as sub- and super-scripts,
respectively, e.g., for $f:A \to B$, we would use $f_A^B$. Given
functions $f$ and $g$ with appropriate domains and co-domains, we
denote their composition with $f \circ g$. We define the canonical
projections $\pi_{A \times B}^A: A \times B \to A$ and
$\pi_{A \times B}^B: A \times B \to B$ by maps $(a,b) \mapsto a$ and
$(a,b) \mapsto b$, respectively. Given a set $A$,
$\identity_A: A \to A$ denotes the identity map.

\section{Koopman Operator Viewpoint of Dynamical
Systems}\label{sec:preliminaries} 
Here, we introduce the Koopman operator viewpoint following~\cite{MB-RM-IM:12} for dynamical systems
when no inputs are present. Consider the discrete-time
system
\begin{align}\label{eq:dynamical-system}
x^+ = T(x), \quad x \in \Xc,
\end{align}
where $T: \Xc \to \Xc$ is a function describing the system's behavior
and $\Xc$ is the state space.  The Koopman operator provides an
alternative viewpoint for describing the behavior of
system~\eqref{eq:dynamical-system} based on examining the evolution of functions defined over the state space. Let $\Fc$ be a vector space (over
field $\cplx$) of complex-valued functions with domain
$\Xc$. Assume $\Fc$ is closed under function composition with
the map $T$. This means that, for all functions
$f: \Xc \to \cplx$ in $\Fc$, one has $f \circ T \in \Fc$. Then, one defines the Koopman operator $\Kc: \Fc \to \Fc$ associated
with~\eqref{eq:dynamical-system} as
\begin{align}\label{eq:Koopman-def}
\Kc f = f \circ T, \quad \forall f \in \Fc.
\end{align}
The action of the Koopman operator $\Kc$ can be viewed as evolving the
value of each function $f$ in $\Fc$, one timestep forward across all
trajectories of~\eqref{eq:dynamical-system}, and encoding the outcome
in the new function $\Kc f \in \Fc$,
\begin{align}\label{eq:Koopman-evolve-trajectories} [\Kc f](x) = f
\circ T(x) = f(x^+), \quad \forall f \in \Fc, \; \forall x \in \Xc.
\end{align}
An important property of the Koopman operator is its linearity, i.e.,
for all $f_1, f_2 \in \Fc$ and $c_1, c_2 \in \cplx$,
\begin{align*}
\Kc (c_1 f_1 + c_2 f_2) = c_1 \Kc f_1 + c_2 \Kc f_2. 
\end{align*}
The linearity of the Koopman operator paves the way to many
interesting applications. We do not deal with them here, but we
mention the most important one: spectral analysis for nonlinear
systems based on a carefully chosen function space $\Fc$ equipped with
additional structure (e.g., norm or inner product). We refer the
reader to e.g.,~\cite{IM:05,IM:21} for more information.

In this paper, our goal is to explore different frameworks to extend
the Koopman operator viewpoint to control systems.

\section{A Cautionary Tale of Naive Koopman Extensions to Control
Systems}\label{sec:naive-Koopman}
In this section we analyze two naive
approaches to build operator viewpoints for control
systems and show their limitations. This sheds light on the
theoretical roadblocks one faces in extending the Koopman operator to control systems.

Consider the control system
\begin{align}\label{eq:control-system}
x^+ = \Tc(x,u), \quad x \in \Xc, \quad u \in \Uc,
\end{align}
where $\Tc: \Xc \times \Uc \to \Xc$ is the function defining the
system's behavior, $\Xc$ is the state space, and $\Uc$ is the input
space.  We do not assume any special form (e.g., control
affine, etc) on system~\eqref{eq:control-system}. We also do not
assume any structure on the sets $\Xc$ and $\Uc$ (not even being
subsets of Euclidean space). As long as $\Tc$
with its given domain and codomain is a well-defined function, for each input, system~\eqref{eq:control-system} has a unique solution starting from each initial condition defined for all time.

\subsection{Naive Approach 1: Simple Composition Operators}\label{sec:naive-1}
One might be tempted to extend the Koopman operator theory to control
systems by defining a composition operator akin
to~\eqref{eq:Koopman-def} in the following sense: let $\Sc_1$ and
$\Sc_2$ be vector spaces over $\cplx$, where $\Sc_1$ is comprised of
complex-valued functions with domain $\Xc$ while $\Sc_2$ is comprised
of complex-valued functions with domain $\Xc \times \Uc$. Moreover,
let $f \circ \Tc \in \Sc_2$ for all $f \in \Sc_1$. Then one can define
a linear composition operator $\Kc^\text{naive} : \Sc_1 \to \Sc_2$
such that
\begin{align}\label{eq:naive-Koopman}
\Kc^\text{naive} f = f \circ \Tc, \quad \forall f \in \Sc_1.
\end{align}
Note that $\Kc^\text{naive}$ is a linear operator and, similarly
to~\eqref{eq:Koopman-evolve-trajectories}, pushes the value of the
function $f$ one step forward in time according to the trajectory
of system~\eqref{eq:control-system}, i.e.,
\begin{align*}
[\Kc^\text{naive} f](x,u) = f (\Tc(x,u)) = f(x^+),
\quad \forall x \in \Xc, \; \forall u \in \Uc .
\end{align*}
Therefore, $\Kc^\text{naive}$ has the two central properties of
Koopman operator: linearity and evolving the function values on system
trajectories. However, it suffers from the ``\emph{curse of domain mismatch}", which invalidates its use for multi-step prediction over the system trajectories: the domain and co-domain of $\Kc^\text{naive}$ are
different. Hence, if $(f: \Xc \to \cplx) \in \Sc_1$, then
$(\Kc^\text{naive} f: \Xc \times \Uc \to \cplx) \in
\Sc_2$. Therefore, one \emph{cannot} apply $\Kc^\text{naive}$ on
$\Kc^\text{naive} f$ to move forward in time for a second timestep,
since the function $\Kc^\text{naive} f$ does not belong to the
domain of $\Kc^\text{naive}$.   

The curse of domain mismatch is a major drawback of the operator
$\Kc^\text{naive}$, since it  does not allow for multi-step prediction.  It is
worth noting that, when $u$ in~\eqref{eq:control-system} represents a
disturbance or process noise rather than a control input, and follows
a probability distribution (which renders the system to be
stochastic), one can resolve the domain mismatch issue by taking the
expectation over $u$, see~\cite{IM-AB:04}\footnote{One can also envision a similar argument for applications in probabilistic reachability analysis.}. However, since here we
focus on open-loop control systems, we do not explore this direction
further and refer the interested reader to~\cite[Section~4]{IM-AB:04}
for a detailed treatment of this approach.

\subsection{Naive Approach 2: Treating the Input as a State}
The second naive approach to extend the Koopman operator theory to
control systems treats the input on the same grounds as the system
state.  One might be tempted to address the curse of domain mismatch by treating the input as an
augmented state and then defining a Koopman operator with
action similar to~\eqref{eq:Koopman-def}. To make this clear, consider
a vector space (over $\cplx$) of functions $\Sc$ comprised of
complex-valued functions with domain $\Xc \times \Uc$ and define a
naive Koopman operator for~\eqref{eq:control-system} denoted by
$\Kc^\text{input-aug}: \Sc \to \Sc$, whose action on each function is
defined similarly to~\eqref{eq:Koopman-evolve-trajectories} as
\begin{align}\label{eq:naive-state-augmented}
&[\Kc^\text{input-aug} f](x,u) = f (x^+, u^+) 
\nonumber
\\
&= f (\Tc(x,u), u^+), \quad \forall f \in \Sc, \; \forall x \in \Xc,
\; 
\forall u, u^+ \in \Uc. 
\end{align}
One already can see an issue with this equation: for an open-loop
system, the input in future timesteps $u^+$ is arbitrary and does not
depend on the state $x$ or the input $u$ at current timestep. For
example, given state and input pair $(x,u) \in \Xc \times \Uc$ at
current time, the state and input pair at the next timestep can be
$(\Tc(x,u), u_1)$ or $(\Tc(x,u), u_2)$ for \emph{different}
$u_1,u_2 \in \Uc$ (if $\Uc$ contains more than one point).  Therefore,
it is not even clear whether the operator $\Kc^\text{input-aug}$ is
well defined. The next result establishes a necessary condition
for $\Kc^\text{input-aug}$ to be well defined, which turns out to be quite restrictive.

\begin{lemma}\longthmtitle{Necessary Condition for Well-defined
$\Kc^\text{input-aug}$}\label{l:well-defined-naive} 
Let the operator $\Kc^\text{input-aug}: \Sc \to \Sc$ associated with
open-loop system~\eqref{eq:control-system} be
well defined. Define the range of the dynamic map $\Tc$
in~\eqref{eq:control-system} as follows
\begin{align}\label{eq:range-map}
\Rc(\Tc) := \setdef{y \in \Xc}{\exists x \in \Xc ,
\exists u \in \Uc \text{ s.t. } y = \Tc(x,u) }. 
\end{align}
Then, at least one of the following hold:
\begin{enumerate}
\item the set $\Uc$ is a singleton;
\item For each $f \in \Sc$, its restriction to $\Rc(\Tc) \times \Uc$ is independent of the second argument.
\end{enumerate}
\end{lemma}
\begin{proof}
At an arbitrary timestep $k$, let $x \in \Xc$ and $u \in \Uc$ be
the state and input of the system. Since operator
$\Kc^\text{input-aug}$ is well defined, it maps every function
$f \in \Sc$ to a \emph{unique} function
$\Kc^\text{input-aug} f \in \Sc$. Moreover, since
system~\eqref{eq:control-system} is open loop, the input at a future
timestep $u^+$ is arbitrary. Therefore, one can choose the
following scenarios: take $u^+ = u_1$ or $u^+ = u_2$ for arbitrary
$u_1, u_2 \in \Uc$. Since equation~\eqref{eq:naive-state-augmented}
holds for all $x \in \Xc$, all $u, u^+ \in \Uc$, and all
$f \in \Sc$, one can write
\begin{align*}
f (\Tc(x,u), u_1) = [\Kc^\text{input-aug} f](x,u) = f (\Tc(x,u), u_2), 
\\
\quad \forall f \in \Sc, \forall x \in \Xc, \forall u, u_1, u_2 \in \Uc.
\end{align*}
Noting that $\Rc(\Tc)$ consists of all $\Tc(x,u)$ for all $(x,u) \in
\Xc \times \Uc$, the previous equation is equivalent to 
\begin{align*}
\restr{f}{\Rc(\Tc) \times \Uc}(y,u_1)
&= \restr{f}{\Rc(\Tc) \times \Uc} (y,u_2),
\\
&\quad \forall f \in \Sc, \, \forall y \in \Rc(\Tc), \, \forall
u_1, u_2 \in \Uc.  
\end{align*}
Since $\restr{f}{\Rc(\Tc) \times \Uc}$ is a function, the equation
above holds under one of the following scenarios:
\begin{enumerate}
\item for all $u_1,u_2 \in \Uc$, $u_1 = u_2$;
\item for all $f \in \Sc$, $\restr{f}{\Rc(\Tc) \times \Uc}$ does not
depend on the second variable (the input at next timestep).
\end{enumerate}
If the former holds, it follows that $\Uc$ contains only one element
and the proof is complete. Otherwise, the latter must hold for all
$f \in \Sc$, concluding the proof.
\end{proof}

Lemma~\ref{l:well-defined-naive} reveals a serious limitation of the
operator defined in~\eqref{eq:naive-state-augmented}. Lemma~\ref{l:well-defined-naive}(a) essentially means that~\eqref{eq:control-system} is not a control system, since $u$ is a constant and can be viewed as a parameter and not as a control
input. Lemma~\ref{l:well-defined-naive}(b) means that $\Kc^\text{input-aug}$ cannot encode multi-step behavior in open-loop systems, as we explain next. 
One should note
that by definition $\Rc(\Tc) \subseteq \Xc$ is forward invariant
(given all possible inputs) and all the trajectories of the system
will end up in $\Rc(\Tc)$ after at most one timestep. Therefore,
the effect of functions $f \in \Sc$ on the system trajectories after
the first timestep can be \emph{completely} captured by their
restriction to $\Rc(\Tc) \times \Uc$. However, based on
Lemma~\ref{l:well-defined-naive}(b), all these restrictions are
independent of the input and \emph{cannot} encode its effect. This
implies the operator $\Kc^\text{input-aug}$
in~\eqref{eq:naive-state-augmented} is unable to capture the
system's behavior for longer than one timestep.

Note the parallelism between this discussion and the curse of domain mismatch for $\Kc^\text{naive}$ in Section~\ref{sec:naive-1}, implying that both operators
$\Kc^\text{naive}$ in~\eqref{eq:naive-Koopman} and
$\Kc^\text{input-aug}$ in~\eqref{eq:naive-state-augmented} fail to
encode multistep trajectories. Next, we remark that despite its
limitation, the operator~\eqref{eq:naive-state-augmented} can be
useful for the analysis of the closed-loop behavior.

\begin{remark}\longthmtitle{Usefulness of $\Kc^\text{input-aug}$ for
Closed-Loop Systems}
Although $\Kc^\text{input-aug}$ in~\eqref{eq:naive-state-augmented} is
of limited utility when it comes to open-loop control systems, this
operator can be useful if the system is closed loop, the input
sequence is fixed in advance, or it complies with a predetermined
dynamics which determines the input \emph{uniquely}. In these cases,
$\Kc^\text{input-aug}$ can be viewed as Koopman operator associated
with the system created by fixing the input structure (which means the
resulting system does not admit a \emph{control} input),
see~\cite{JLP-SLB-JNK:18}. \oprocend
\end{remark}
\smallskip

The naive
approaches~\eqref{eq:naive-Koopman}-\eqref{eq:naive-state-augmented}
to extend the Koopman operator theory to open-loop control systems
reveal a major difficulty arising from the fact that the input has a
fundamentally different role compared to the state of a control
system: state follows a prescribed dynamic map while input in an open
loop system is arbitrary and can change the behavior of the dynamic
map itself.

\begin{remark}\longthmtitle{What Went Wrong? Input $\neq$
State}\label{r:input-not-state}
The difference between the roles of state and input highlights a
subtle but crucial point in extending Koopman operator theory to
control. The Koopman operator in~\eqref{eq:Koopman-def} is simply an
alternative representation of the dynamical map $T$,
cf.~\eqref{eq:dynamical-system}, in an appropriate function
space. Without the underlying map $T$, the Koopman operator has no
meaning. Unlike the state, which evolves according to a map, the
input in an open-loop control system does \textbf{not} follow any
predefined evolution. Hence, there is no map that can be directly
represented as a linear operator. Consequently, one must develop
other approaches to properly capture the effect of input, which is
what we discuss next.  \oprocend
\end{remark}

\section{Formal Extensions of Koopman Operator Theory to Control
Systems}\label{sec:formal-extensions}
The difficulty of extending the Koopman operator theory to control
systems lies in the fact that input to an open-loop system is
arbitrary. Changing the input sequences may drastically alter the
system behavior. Therefore, in extending the Koopman operator theory
to control systems, one should take into account all possible
behaviors arising from different input sequences.

To the best of our knowledge, there are two general
extensions\footnote{By a general extension, we mean extensions based
on general operator-theoretic descriptions for general (not
necessarily control affine) nonlinear systems. There exist methods
based on finite-dimensional lifting approaches (e.g.,
super-linearization~\cite{JLP-SLB-JNK:16} or
bilinearization~\cite{DG-DAP:21}) and operator theoretic methods
based on specific assumptions on dynamic maps (e.g., control
affine~\cite{SP-SEO-CWR:20}), which we do not discuss here.}  of the
Koopman operator theory to control (not necessarily control-affine)
systems and both turn the control system into systems without
input. The first approach~\cite{MK-IM-automatica:18} achieves this by
considering all possible infinite input sequences, while the second
approach~\cite{MH-JC:25-auto} considers all possible systems that one
can build by setting the input in~\eqref{eq:control-system} to be a
constant. We describe both next.

\subsection{Koopman Operator via Infinite Input
Sequences}\label{sec:Koopman-infinite-sequence}
Following~\cite{MK-IM-automatica:18}, consider the space $\ell (\Uc)$
comprised of all infinite sequences
$\mathbf{u}:= (u_n)_{n = 0}^\infty$, where $u_n \in \Uc$ for all
$n \in \naturals_0$. Then, one can define a system without input on
the set $\Xc \times \ell(\Uc)$ as
\begin{align}\label{eq:system-infinite-sequence-large}
\begin{bmatrix}
x
\\
\mathbf{u}
\end{bmatrix}^+
=
\begin{bmatrix}
\Tc(x ,
\mathbf{u}(0))
\\
S_{\text{left}}
\mathbf{u}
\end{bmatrix}, 
\end{align}
where $S_{\text{left}} : \ell(\Uc) \to \ell(\Uc)$ is the left shift
operator defined by the mapping
$\big(u(0) , u(1), \ldots \big) \mapsto \big(u(1) , u(2), \ldots
\big)$.

For convenience, we denote the system
in~\eqref{eq:system-infinite-sequence-large} by a tuple notation
\begin{align}\label{eq:system-infinite-sequence}
(x, \mathbf{u})^+ = \Tc^{\infty} (x, \mathbf{u}) :=
(\Tc(x , \mathbf{u}(0)), S_{\text{left}} \mathbf{u}), 
\nonumber
\\ 
\forall (x, \mathbf{u}) \in \Xc \times \ell(\Uc).
\end{align}
Note that the system defined by
$\Tc^{\infty} : \Xc \times \ell (\Uc) \to \Xc \times \ell (\Uc)$ is
now a system without input and admits a Koopman operator similarly
to~\eqref{eq:Koopman-def}. Let $\Fc^{\infty}$ be a vector space (over
$\cplx$) of complex-valued functions with domain
$\Xc \times \ell (\Uc)$ that is closed under composition with
$\Tc^{\infty}$. Then, we define the following Koopman operator
$\Kc^\infty : \Fc^\infty \to \Fc^\infty$ as
\begin{align}\label{eq:Koopman-infinite-sequence}
\Kc^\infty f = f \circ \Tc^{\infty}, \quad \forall f \in \Fc^\infty.
\end{align}
Unlike the approaches in Section~\ref{sec:naive-Koopman}, the operator
$\Kc^{\infty}$ is always well-defined as long as the function space
$\Fc^\infty$ is closed under composition with $\Tc_{\infty}$.  Note
that the operator $\Kc^\infty$ is the Koopman operator associated with
the system~\eqref{eq:system-infinite-sequence}, which is not the
control system~\eqref{eq:control-system}. Therefore, we need to show
that it can capture the information about the control
system~\eqref{eq:control-system}. To achieve this goal we rely on the
following notion.

\begin{definition}\longthmtitle{Control-Independent Functions in
$\Finf$ and State
Component}\label{def:control-independent-func-inf} 
The function $f: \Xc \times \ell(\Uc) \to \cplx$ in $\Finf$ is
\emph{control-independent} if
$f(x,\mathbf{u_1}) = f(x,\mathbf{u_2})$ for all $x \in \Xc$ and all
$\mathbf{u_1}, \mathbf{u_2} \in \ell(\Uc)$. Alternatively, $f$ can
be decomposed as
\begin{align*}
f(x,\mathbf{u}) = f_\Xc(x) 1_{\ell(\Uc)}(\mathbf{u}), \quad
\forall (x,\mathbf{u}) \in \Xc \times \ell(\Uc),
\end{align*}
where $1_{\ell(\Uc)}: \ell(\Uc) \to \cplx$ is the constant function
equal to 1. We call $f_\Xc: \Xc \to \cplx$ the \emph{state
component} of $f$. We denote by $\Finfci$ the set of all control-independent
functions in $\Finf$.  \oprocend
\end{definition}
\smallskip

Although the domain of control-independent functions is
$\Xc \times \ell(\Uc)$, they only depend on the part corresponding to
the state of the original control system~\eqref{eq:control-system}
(which is $\Xc$) and discard the information from
$\ell(\Uc)$. Therefore, control-independent functions can be
completely captured via their state component.

The next result shows that control-independent functions recover the
information of the trajectories of original control
system~\eqref{eq:control-system} through their state components.

\begin{lemma}\longthmtitle{Encoding Information of Control
System~\eqref{eq:control-system} via the Action of $\Kc^\infty$ on
Control-Independent
Functions}\label{l:capture-system-control-ind-inf} 
Let $f \in \Finfci$ and let $\{x_k\}_{k \in \naturals_0}$ be the
trajectory of~\eqref{eq:control-system} from initial condition $x_0$
with input sequence $\mathbf{u} = (u_0, u_1, \ldots)$. Then,
\begin{align*}
[(\Kc^\infty)^k f](x_0,\mathbf{u} ) = f_\Xc(x_k),
\quad \forall k \in \naturals_0,
\end{align*}
where $(\Kc^\infty)^k$ is the composition,  $k$ times, of $\Kc^\infty$
with itself.  \oprocend
\end{lemma}
\smallskip

The proof of Lemma~\ref{l:capture-system-control-ind-inf}, omitted for
space reasons, trivially follows from an inductive process on
$\Kc^\infty$ and applying the decomposition in
Definition~\ref{def:control-independent-func-inf}.
Lemma~\ref{l:capture-system-control-ind-inf} shows that one can use
control-independent functions in conjunction with the operator
$\Kc^\infty$ to encode information from the trajectories
of~\eqref{eq:control-system}. A simple example of control-independent
functions are state observables: given $\Xc \subseteq \real^n$, let
$\Finf$ contain the function $h_i$, $h_i(x, \mathbf{u}) = x^{(i)}$,
where $x^{(i)}$ is the $i$th element of the state
of~\eqref{eq:control-system}, $i \in \until{n}$. Then, one can use
Lemma~\ref{l:capture-system-control-ind-inf} to extract the $i$th
component of the control system's trajectories for all time using the
operator $\Kc^\infty$.

The state space of system~\eqref{eq:system-infinite-sequence} is
$\Xc \times \ell(\Uc)$; therefore, to evaluate the dynamics, one
requires a point in $\Xc$, the state space of~\eqref{eq:control-system}, and an \emph{infinite} sequence of
inputs. Given that the domain of functions in $\Finf$ is also
$\Xc \times \ell(\Uc)$, the same requirements  apply to evaluate them. Hence, it is not possible to truncate the
input sequence for function evaluation. The infinite-input sequences
are the reason why the operator $\Kc^{\infty}$ can encode complete
information about system~\eqref{eq:control-system} and does not suffer
from the same issues as the approaches in
Section~\ref{sec:naive-Koopman}.

\begin{remark}\longthmtitle{Finite-Dimensional Representations for
Infinite-Input  Sequences Framework}\label{rem:finite-infinite}
Applying Koopman-based techniques on digital computers usually
involves approximations on finite-dimensional spaces often achieved
through projections, which lead to information loss. Therefore, for
such finite-dimensional forms, it is critical to provide error
bounds on the prediction of all functions in the subspace and derive
methods to learn such models from data.  In the case of control
systems, deriving finite-dimensional forms is considerably
involved. Many works draw inspiration from the idea of
\emph{lifting} to obtain such forms. Lifted linear and bilinear
models are common choices, motivated by their simplicity and ease of
use. However, the connection of these models to the Koopman-based
description is not always clear, and in some cases may be absent
altogether\footnote{It is worth noting that for
\emph{continuous-time control-affine} systems, the lifted bilinear
form (while still an approximation) can be used effectively for
finite-horizon prediction in model predictive control schemes and
related feedback designs. For detailed error bounds, and
closed-loop guarantees, we refer the reader
to~\cite{RS-KW-IM-JB-MS-FG:25} and references therein.}.  Given
that the functions in $\Finf$ fuse the information of states and
infinite input sequences together, directly finding
finite-dimensional forms for this framework is more difficult and,
to the best of our knowledge, has not been done in the literature
for general (not necessarily control-affine) control systems.
\oprocend
\end{remark}

\subsection{Koopman Control Family (KCF)}\label{sec:KCF}
The second approach also turns the control
system~\eqref{eq:control-system} into systems without input, so that
Koopman operators akin to~\eqref{eq:Koopman-def} can be
used. Following~\cite{MH-JC:25-auto}, this construction is done by
considering all possible systems that one can build by fixing the
input to be constant in~\eqref{eq:control-system}. This leads to
\begin{align}\label{eq:constant-input-family}
x^+ = \Tc_{u^*} (x) := \Tc(x, u \equiv u^*), \quad u^* \in \Uc,
\end{align}
where the set $\{\Tc_{u^*} : \Xc \to \Xc \}_{u^* \in \Uc}$ forms a
family of systems without input. Note that the domain and co-domain of
this parametric family match, since they are no longer control
systems. Therefore, we can assign to each map a Koopman operator
similarly to~\eqref{eq:Koopman-def}. Let $\Fc$ be a vector space of
functions (over $\cplx$) of complex-valued functions that are closed
under composition with the members of family
$\{\Tc_{u^*} : \Xc \to \Xc \}_{u^* \in \Uc}$. Then, we define the
Koopman control family (KCF) as
$\{\Kc_{u^*} : \Fc \to \Fc\}_{u^* \in \Uc}$, where
\begin{align}\label{eq:Koopman-control-family}
\Kc_{u^*} f = f \circ \Tc_{u^*}, \quad \forall f \in \Fc, \quad
\forall u^* \in \Uc. 
\end{align}

The following result shows how the KCF encodes the trajectories
of~\eqref{eq:control-system}.

\begin{lemma}\longthmtitle{Encoding Information of Control System~\eqref{eq:control-system} via
KCF}\label{l:encoding-trajectories-KCF}
Let $g \in \Fc$ and let $\{x_k\}_{k \in \naturals_0}$ be the
trajectory of~\eqref{eq:control-system} from initial condition $x_0$
with input sequence $\mathbf{u} = (u_0, u_1, \ldots)$. Then,
\begin{align*}
[\Kc_{u_0} \Kc_{u_1} \ldots \Kc_{u_{k-1}} g](x_0) = g(x_k) , \quad
\forall k \in \naturals_0. 
\eqoprocend 
\end{align*}
\end{lemma}

The proof of Lemma~\ref{l:encoding-trajectories-KCF} trivially follows
from the definition of KCF and is omitted. One can use
Lemma~\ref{l:encoding-trajectories-KCF} to fully extract the system
trajectories via state observables if they belong to~$\Fc$: given
$\Xc \subseteq \real^n$, let $\Fc$ contain the function $o_i$,
$o_i(x) = x^{(i)}$, where $x^{(i)}$ is the $i$th component of $x$,
$i \in \until{n}$. Then, one can use
Lemma~\ref{l:encoding-trajectories-KCF} to fully recover the $i$th
component of the trajectory of system~\eqref{eq:control-system} for
all time.

\begin{remark}\longthmtitle{KCF and Switch-Based Koopman Control}
Although not as formal extensions of Koopman operator theory to
control systems, the idea of controlling nonlinear systems by
considering finitely many possible constant inputs and switching
between finite-dimensional approximations of the associated Koopman
operators was already explored in the literature (see,
e.g.,~\cite{SP-SK:19}) prior to the introduction of
KCF~\cite{MH-JC:25-auto}.  Moreover, the
works~\cite{AS-DE:17,AS-AM-DE:18} also considered piece-wise
constant input control for special types of systems.  One can view
KCF as a formal generalization of such approaches to the case of
abstract vector spaces, uncountable input sets, and arbitrary
control systems.  We refer the reader to~\cite{MH-JC:25-auto} for a
detailed discussion of differences and how the finite-dimensional
form derived from KCF generalizes lifted linear, lifted bilinear,
and lifted switched-linear forms in the literature.  \oprocend
\end{remark}

\begin{remark}\longthmtitle{Finite-Dimensional Representations for KCF}\label{rem:finite-KCF}
Given a finite-dimensional subspace $\Sc \in \Fc$, let
$\Psi: \Xc \to \cplx^{\dim(\Sc)}$ be a vector-valued function whose
elements span $\Sc$. Then the finite-dimensional from for KCF,
termed ``input-state separable'' model, is as follows
\begin{align}\label{eq:input-state-separable}
\Psi (x^+) = \Psi \circ \Tc(x,u) \approx \Ac(u) \Psi(x), \quad
\forall x \in \Xc, u \in \Uc, 
\end{align}
where $\Ac: \Uc \to \cplx^{{\dim(\Sc) \times \dim(\Sc)}}$ is a
matrix-valued function. The form in~\eqref{eq:input-state-separable}
is a result of a tight (necessary and sufficient) condition which
cannot be relaxed~\cite[Theorem~4.3]{MH-JC:25-auto}.  The
input-state separable model is linear in $\Psi(x)$ (often referred
to as lifted state) but generally nonlinear in input $u$. The reason
for this is as follows: Koopman-based methods represent nonlinear
dynamic maps via linear operators; however, unlike the state $x$,
which is governed by the dynamic map~\eqref{eq:control-system}, the
free input $u$ is arbitrary an is not governed by a map; therefore,
one cannot represent the evolution of input via a linear
operator. This distinction is the reason why the naive approaches in
Section~\ref{sec:naive-Koopman} fail,
cf.~Remark~\ref{r:input-not-state}. Interestingly, the commonly used
lifted linear, bilinear, and switched linear models mentioned above
are all special cases of the input-state separable form.
We refer the interested reader to~\cite{MH-JC:25-auto} for the
detailed derivation of input-state separable forms, error bounds for
function prediction using the notion of invariance proximity, as
well as data-driven learning with accuracy guarantees.  \oprocend
\end{remark}

\section{Motivation and Problem Statement}\label{sec:problem-statement}
The question we address in this paper is how the frameworks described in 
Sections~\ref{sec:Koopman-infinite-sequence} and~\ref{sec:KCF} are related to each other. Based on our discussion, one would expect them to be equivalent. However, this is not always the case and the
answer depends on the choice of function spaces, as the following simple example shows.

\begin{example}\longthmtitle{Dependence on Function spaces}\label{ex:simple}
Let $\Finf = \Span(1_{\Xc \times \ell(\Uc)})$ be the space of
constant functions on $\Xc \times \ell(\Uc)$ and $\Fc = B(\Xc)$ be
the space of bounded functions on $\Xc$ (both on the field $\cplx$). One
can readily check that the operators in both frameworks are well
defined. However, in this case, the operator $\Kc^\infty$ is trivial
and does not capture any information about the system since constant
functions do not change under composition with any map. On the other hand, $\Fc$ is richer and can
capture some information about the system dynamics. One could switch
the construction of spaces to $\Finf = B(\Xc \times \ell(\Uc))$ and
$\Fc = \Span(1_\Xc)$, in which case, KCF will not capture any
information, while $\Kc^\infty$ captures some.
\oprocend
\end{example}

Example~\ref{ex:simple} makes it clear that the function spaces $\Finf$ and
$\Fc$ should satisfy certain conditions for the frameworks to capture complete information about the system, as a necessary step to establish their equivalence.

\begin{problem}\longthmtitle{Equivalence of Formal Extensions of
Koopman Theory to Control Systems} 
Given the frameworks of the Koopman operator via infinite input
sequences ($\Kc^\infty$) and the Koopman control family (KCF):
\begin{enumerate}
\item provide conditions on the function spaces $\Fc$ and $\Fc^\infty$
such that $\Kc^\infty$ and KCF can be converted to each other only
via compositions with linear operators between the functions spaces;
\item given access to $\Kc^\infty$ or KCF, provide constructive
recipes to predict the action of the operators on functions in the
other framework. This task should be performed at the highest level
of generality with no reliance on specific structure of function
spaces (e.g., topology, metric, norm, inner product, etc) or
specific structures on the dynamics (e.g., being control affine).
\oprocend
\end{enumerate}
\end{problem}

Beyond the mathematical relevance of establishing the equivalence between the two different frameworks, there are also important practical implications. Each framework handles system information in fundamentally different
ways. Depending on the application, one framework may therefore be
more convenient than the other.  In practice, function spaces are
often equipped with additional structures --such as norms or inner
products-- or are subject to constraints, for example through the use
of tensor product spaces. For a given application, one can select the
framework best suited to the task at hand. This choice is dictated
directly by the nature of the application.  For applications concerned
with the general study of control systems under all possible inputs --
such as the analysis of invariant quantities, the characterization of
unreachable sets, or theoretical investigations into the spectral
properties of operators -- the infinite-input sequence framework is a
natural choice, since it relies on a single operator. By contrast, in
applications involving finite-dimensional representations or
finite-time trajectories, such as finite-horizon prediction in model predictive control (MPC) or data-driven learning via trajectory data,
the
fact that the KCF framework avoids the need to work with
infinite-input sequences when evaluating functions in the function
space together with the availability of the input-state separable
form~\eqref{eq:input-state-separable} makes it more attractive.

\section{Connections Between Dynamical Systems in Different Frameworks}
Both approaches in Section~\ref{sec:formal-extensions} use infinitely
many objects to provide Koopman-based descriptions: one 
relies on infinite input sequences while the other utilizes
infinitely many constant-input systems. This infinite cardinality
leads to practical difficulties. To address this, here we introduce a different, easier-to-work-with system and 
show it indirectly captures all the
information needed to reconstruct the original system trajectories.
We then utilize this system to study and unify the aforementioned
extensions of Koopman operator theory to control systems.

\subsection{Augmented System and its Associated Koopman Operator}\label{sec:kaug}
We utilize the control system~\eqref{eq:control-system} to synthesize
a new system without input, cf.~\cite{MH-JC:25-auto},
\begin{align*}
\begin{bmatrix}
x \\ u
\end{bmatrix}^+
=
\begin{bmatrix}
\Tc(x , u)
\\
u
\end{bmatrix},
\quad x \in \Xc, \quad
u \in \Uc. 
\end{align*}
Note that in this system, $u$ is part of the state and does not
evolve in time. For convenience, we use the following notation for the
augmented system
\begin{align}\label{eq:augmented-system}
(x,u)^+ = (\Tc(x,u), u) =: \Taug (x,u),
\end{align}
where $\Taug : \Xc \times \Uc \to \Xc \times \Uc$ is the function
defining the augmented system. Note that~\eqref{eq:augmented-system}
does not have a control input; hence, given a suitable function space,
admits a well-defined Koopman operator similarly
to~\eqref{eq:Koopman-def}. Let $\Faug$ be a vector space (over
$\cplx$) of complex-valued functions with domain $\Xc \times \Uc$ that
is closed under composition with $\Taug$. Then, we define the
augmented Koopman operator $\Kaug : \Faug \to \Faug$ as
\begin{align}\label{eq:augmented-Koopman}
\Kaug f = f \circ \Taug, \quad \forall f \in \Faug.
\end{align}
Unlike the system with infinite input
sequences~\eqref{eq:system-infinite-sequence} or the family of
constant input systems~\eqref{eq:constant-input-family}, the augmented
system~\eqref{eq:augmented-system} does not rely on any infinite
objects and therefore is easier to study and evaluate; hence, it
serves as an effective intermediary between the two frameworks.

Similarly to Definition~\ref{def:control-independent-func-inf}, we
define a set of control-independent functions in $\Faug$ to connect
the action of $\Kaug$ to the trajectories of control
system~\eqref{eq:control-system}.

\begin{definition}\longthmtitle{Control-Independent Functions in
$\Faug$ and State
Component}\label{def:control-independent-func-aug}
The function $f: \Xc \times \Uc \to \cplx$ in $\Faug$ is
\emph{control-independent} if $f(x,u_1) = f(x,u_2)$ for all
$x \in \Xc$ and all $u_1, u_2 \in \Uc$. Alternatively, $f$ can be
decomposed as
\begin{align*}
f(x,u) = f_\Xc(x) 1_{\Uc}(u), \quad \forall (x,u) \in \Xc \times \Uc,
\end{align*}
where $1_{\Uc}: \Uc \to \cplx$ is a constant function equal to
1. We call $f_\Xc: \Xc \to \cplx$ the \emph{state component}
of $f$. We denote the set of all control-independent functions
in $\Faug$ by $\Faugci$.  \oprocend
\end{definition}
\smallskip

The next result shows that one can use the control-independent functions in
conjunction with $\Kaug$ to extract information about  the
trajectories of~\eqref{eq:control-system}.

\begin{lemma}\longthmtitle{Encoding Single-Step Information of Control
System~\eqref{eq:control-system} via the Action of $\Kaug$ on
Control-Independent
Functions}\label{l:capture-system-control-ind-aug}
For $f \in \Faugci$, we have
\begin{align*} [\Kaug f](x,u) = f_\Xc(x^+), \quad \forall (x,u) \in
\Xc \times \Uc. \eqoprocend
\end{align*}
\end{lemma}

The proof of Lemma~\ref{l:capture-system-control-ind-aug} directly
follows from the definitions. It is crucial to note that unlike
Lemma~\ref{l:capture-system-control-ind-inf}, the prediction in
Lemma~\ref{l:capture-system-control-ind-aug} only holds for one time
step. Hence, $\Kaug$ is not a Koopman operator associated with control
system~\eqref{eq:control-system}. In fact, one can easily show that,
in general, $\Kaug$ cannot directly capture the action
of~\eqref{eq:control-system} for longer than a single timestep and
suffers from similar issues as the naive approaches in
Section~\ref{sec:naive-Koopman}. However, $\Kaug$ has several important properties in
\emph{indirectly} parameterizing the KCF and deriving
finite-dimensional forms for Koopman-based approaches, cf.~\cite{MH-JC:25-auto}.

\subsection{All Introduced Input-Free Systems Capture the Behavior of
the Control System}

The next result shows that one can completely, but indirectly, recover
the original system's behavior from the input-free dynamical systems
introduced so far.

\begin{proposition}\longthmtitle{Recovering the Original System's
Behavior from Different Input-Free Systems}\label{p:recoving-trajectory}
Let $\{x_k\}_{k \in \naturals_0}$ be the trajectory of the original
system~\eqref{eq:control-system} with initial condition $x_0$ for the
input sequence $\mathbf{u} = (u_0, u_1, \ldots)$.  Then, the
trajectory can be fully recovered from the system with infinite
input-sequences~\eqref{eq:system-infinite-sequence}, the family of
constant-input systems~\eqref{eq:constant-input-family}, and the
augmented system~\eqref{eq:augmented-system} as follows: for all
$k \in \naturals$,
\begin{enumerate}
\item \label{eq:prop-infinite-sequence}
$x_k = [\pi^{\Xc}_{\Xc \times \ell(\Uc)} \circ
\big(\Tc_{\infty}\big)^{k}] (x_0, \mathbf{u})$;
\item \label{eq:prop-constant-input}
$x_k = \Tc_{u_k} \circ \cdots \circ \Tc_{u_0} (x_0)$;
\item \label{eq:prop-augmented}
$x_k = \pi^{\Xc}_{\Xc \times \Uc} \circ \Taug(x_{k-1}, u_{k-1})$. \oprocend
\end{enumerate}
\end{proposition}
\smallskip

The proof of Proposition~\ref{p:recoving-trajectory} can be done via
direct calculation. One should note the difference in
Proposition~\ref{p:recoving-trajectory}(c) with the other parts, which
is rooted in the fact that $\Taug$ encodes the information of the
system~\eqref{eq:control-system} only for one timestep.
Proposition~\ref{p:recoving-trajectory} shows that all the introduced
frameworks can indirectly capture complete information about the
original system. Therefore, by properly choosing the function spaces,
one should be able to connect their associated Koopman-based
structures. This is what we tackle in the next sections.

\section{Connecting the Function Spaces}
Here, we take a step towards connecting the different
operator-theoretic descriptions by first connecting their associated function
spaces.  As summarized in Table~\ref{tab:framework-comparison}, the
domains of the each set of functions are the state
space of the corresponding dynamics and are therefore different. 
We then need operations to change the
domain of functions and move between function spaces. To achieve this
goal, we use $\Faug$ as an intermediary to connect $\Fc^\infty$
and~$\Fc$.

\begin{table}[h]
\centering
\renewcommand{\arraystretch}{1.5}
\begin{tabular}{|c|c|c|c|}
\hline
\textbf{Framework} & \textbf{Dynamics} & \textbf{State Space} & \textbf{Function Space} \\
\hline
\text{Infinite Sequences} & $\Tc^\infty$ & $\Xc \times \ell(\Uc)$ & $\Fc^\infty$ \\
\hline
\text{KCF} & $\{\Tc_{u^*}\}_{u^* \in \Uc}$ & $\Xc$ & $\Fc$ \\
\hline
\text{Augmented} & $\Taug$  & $\Xc \times \Uc$ & $\Faug$ \\
\hline
\end{tabular}
\caption{Comparison of dynamics, state space, and function spaces in
different frameworks.}
\label{tab:framework-comparison}
\end{table}

\subsection{Connecting $\Fc^\infty$ and $\Faug$}\label{s:functions-inf-aug}
To connect $\Fc^\infty$ and $\Faug$, we provide two operations to
switch the domain of functions between $\Xc \times \ell(\Uc)$ and
$\Xc \times \Uc$.

\begin{definition}\longthmtitle{Domain Restriction and Extension between $\Xc \times \ell(\Uc)$ and
$\Xc \times \Uc$}\label{def:infinite-augmented-res-aug}
Let $f: \Xc \times \ell(\Uc) \to \cplx$ be a function in $\Fc^\infty$ and $g: \Xc \times \Uc \to \cplx$ a function in $\Fc$. Then,
\begin{enumerate}
\item the \emph{domain restriction} of $f$ to $\Xc \times \Uc$, denoted $\restr{f}{\Xc \times \Uc}$, is
\begin{align*}
\restr{f}{\Xc \times \Uc} (x,u) := f (x, (u, u, \ldots)), \quad
\forall (x,u) \in \Xc \times \Uc. 
\end{align*}

\item the \emph{domain extension} of $g$ to $\Xc \times \ell(\Uc)$, denoted
$g^\infty$, is
\begin{align*}
g^\infty (x, \mathbf{u}) := g(x, \mathbf{u}(0)), \quad \forall (x,
\mathbf{u}) \in \Xc \times \ell(\Uc). 
\eqoprocend
\end{align*}
\end{enumerate}
\end{definition}
\smallskip

Definition~\ref{def:infinite-augmented-res-aug} provides an intuitive way to connect the function spaces. The following provides an equivalent description based on mappings, which later will be more convenient to reason with the operator-theoretic
descriptions. 

\begin{lemma}\longthmtitle{Mapping-based Connection between $\Xc \times \ell(\Uc)$ 
and $\Xc \times \Uc$}\label{l:mapping-based-infinite-augmented}
Define the mappings\footnote{The letters $R$ and $E$ in the name of
mappings~\eqref{eq:mapping-based-infinite-augmented} refer to
restriction and extension, respectively. We use similar
terminology for maps and operators throughout the paper.}
$\Rinftoaug : \Xc \times \ell(\Uc) \to \Xc \times \Uc$ and
$\Eaugtoinf : \Xc \times \Uc \to \Xc \times \ell(\Uc)$ as
\begin{subequations}\label{eq:mapping-based-infinite-augmented}
\begin{align}
\Rinftoaug (x, \mathbf{u}) = (x, \mathbf{u}(0)), \quad \forall
(x, \mathbf{u}) \in \Xc \times \ell(\Uc), 
\\
\Eaugtoinf (x, u) = (x, (u, u, \ldots)), \quad \forall (x, u)
\in \Xc \times \Uc. 
\end{align}
\end{subequations}
Then,
\begin{enumerate}
\item $\restr{f}{\Xc \times \Uc} = f \circ \Eaugtoinf $, for all
$f \in \Fc^\infty$;
\item  $g^\infty = g \circ \Rinftoaug $, for all $g \in \Faug$;
\item $\Rinftoaug \circ \Eaugtoinf = \identity_{\Xc \times \Uc}$.
\oprocend
\end{enumerate}
\end{lemma}

We omit the  proof of Lemma~\ref{l:mapping-based-infinite-augmented}, which
follows from the definitions. This result allows us to move
between the functions spaces $\Fc^\infty$ and $\Faug$ via linear
composition operators, as we explain next.

\begin{proposition}\longthmtitle{Linear Operator Connection Between
$\Fc^\infty$ and
$\Faug$}\label{p:linear-connection-infinite-augmented}
Assume $\Faug$ and $\Fc^\infty$ satisfy:
\begin{enumerate}[label=($\mathtt{C}$\roman*)]
\item $\restr{f}{\Xc \times
\Uc}  \in \Faug $, for all ${f \in \Finf}$;
\item $g^\infty \in \Finf$, for all $g \in \Faug$.
\end{enumerate}
Define the operators $\RopFinftoFaug : \Finf \to \Faug$
and $\EopFaugtoFinf : \Faug \to \Finf$ as
\begin{align*}
\RopFinftoFaug f = f \circ \Eaugtoinf, \qquad 
\EopFaugtoFinf g = g \circ \Rinftoaug .
\end{align*}
Then,
\begin{enumerate}
\item $\RopFinftoFaug$ and $\EopFaugtoFinf$ are well-defined and linear;
\item $\RopFinftoFaug f = \restr{f}{\Xc \times \Uc}$, for all $f \in \Finf$;
\item $\EopFaugtoFinf g = g^\infty$, for all $g \in \Faug$;
\item $\RopFinftoFaug \EopFaugtoFinf = \identity_{\Faug}$;
\item $\EopFaugtoFinf(\Faugci) = \Finfci$;
\item $\RopFinftoFaug(\Finfci) =\Faugci$.
\end{enumerate}
\end{proposition}
\begin{proof}
(a) For $f \in \Finf$, using
Lemma~\ref{l:mapping-based-infinite-augmented}, we have
$ \RopFinftoFaug f = \restr{f}{\Xc \times
\Uc}$. Therefore, $\RopFinftoFaug f \in \Faug$ and the operator
is well-defined. A similar argument holds for $\EopFaugtoFinf$. Linearity of the operators directly
follows from their definition.

(b)-(c) The proof follows from
Lemma~\ref{l:mapping-based-infinite-augmented}.

(d) For $g \in \Faug$, using
Lemma~\ref{l:mapping-based-infinite-augmented}(c), one can write
$\RopFinftoFaug \EopFaugtoFinf g = g \circ \Rinftoaug \circ
\Eaugtoinf = g \circ \identity_{\Xc \times \Uc} = g$.

(e) We first prove $\EopFaugtoFinf(\Faugci) \subseteq \Finfci$. Let
$g \in \Faugci$. Then, by
Definition~\ref{def:control-independent-func-aug},
$g(x,u_1) = g(x,u_2)$, for all $x \in \Xc$, $u_1, u_2 \in \Uc$.
Hence, we can write
$[\EopFaugtoFinf g](x,\mathbf{u_1}) = g(x,\mathbf{u_1}(0)) =
g(x,\mathbf{u_2}(0)) = [\EopFaugtoFinf g](x,\mathbf{u_2})$ for all
$x \in \Xc$ and all $\mathbf{u_1}, \mathbf{u_2} \in
\ell(\Uc)$. Hence $\EopFaugtoFinf g \in \Finfci$.

Next, we prove $\Finfci \subseteq \EopFaugtoFinf(\Faugci)$. Given
$f \in \Finfci$,
let
$g_f = \restr{f}{\Xc \times \Uc}$. Since $f$ is
control-independent, then $g_f = \restr{f}{\Xc \times \Uc}$ is also
control-independent. Hence, we only need to show
$\EopFaugtoFinf g_f = f$. For
$(x, \mathbf{u}) \in \Xc \times \ell(\Uc)$,
\begin{align*}
[\EopFaugtoFinf g_f](x, \mathbf{u}) =
\restr{f}{\Xc \times \Uc}(x, \mathbf{u}(0))
&= f(x,
(\mathbf{u}(0), \mathbf{u}(0), \ldots))
\\
&=f(x,\mathbf{u}),
\end{align*}
where
in the last equality we use that $f$ is
control-independent. This proves
$\Finfci \subseteq \EopFaugtoFinf(\Faugci)$, and we conclude $\Finfci = \EopFaugtoFinf(\Faugci)$.

(f) Directly follows from (d) and (e).
\end{proof}

Proposition~\ref{p:linear-connection-infinite-augmented} provide
conditions under which we can move between the function spaces
$\Fc^\infty$ and $\Faug$ via linear operators. These conditions on
function spaces will later be essential in establishing an equivalence
between the frameworks.  We note that the
restriction operator $\RopFinftoFaug$ is built with the extension
map $\Eaugtoinf$ and the extension operator $\EopFaugtoFinf$ is built
with the restriction map $\Rinftoaug$. This is because the operators
act from the left and the maps are composed from the right.

\subsection{Connecting $\Fc$ and $\Faug$}\label{s:functions-ori-aug}
Here, we connect the
function spaces $\Fc$ and $\Faug$. This connection is more complicated than connecting $\Fc^\infty$
and $\Faug$ because the KCF is generally comprised of uncountably many
operators. In our previous work~\cite{MH-JC:25-auto}, we have analyzed the connection between KCF and the augmented operator, but the perspective here is slightly different, with milder conditions suited for the
problem at hand. We start by providing suitable
notions of domain restriction and extension between $\Xc \times \Uc$
and~$\Xc$.

\begin{definition}\longthmtitle{Domain Restriction and Extension between $\Xc \times \Uc$ and
$\Xc$}\label{def:augmented-original-res-aug}
Let $f: \Xc \times \Uc \to \cplx$ be a function in $\Faug$ and $g: \Xc \to \cplx$ a function in $\Fc$. Then,
\begin{enumerate}
\item the family of \emph{constant-input domain restrictions} of $f$ to $\Xc$, denoted 
$\{\restr{f}{x, u \equiv u^*} \}_{u^* \in \Uc}$, is, for each $u^* \in \Uc$,
\begin{align*}
\restr{f}{\Xc, u \equiv u^*} (x) := f (x, u^*), \quad \forall x
\in \Xc . 
\end{align*}

\item the \emph{input-independent domain extension} of $g$ to $\Xc \times \Uc$,
denoted $g_e$, is
\begin{align*}
g_e (x, u) := g(x), \quad \forall (x, u) \in \Xc \times \Uc.
\eqoprocend
\end{align*}
\end{enumerate}
\end{definition}
\smallskip

We note that there is a family of restrictions in Definition~\ref{def:augmented-original-res-aug}(a), instead of the one restriction in 
Definition~\ref{def:infinite-augmented-res-aug}(a).
The next result, analogous to Lemma~\ref{l:mapping-based-infinite-augmented}, provides mapping-based
descriptions of the domain restriction and extension.

\begin{lemma}\longthmtitle{Mapping-based Connection Between
$\Xc \times \Uc$ and
$\Xc$}\label{l:mapping-based-augmented-original}
Define the mapping $\Raugtoorig : \Xc \times \Uc \to \Xc$
and the family of mappings
$\{\Eorigtoaug{u^*} : \Xc \to \Xc \times \Uc \}_{u^* \in \Uc}$ as
\begin{subequations}\label{eq:maps-e-r-aug-orig}
\begin{align}
&\Raugtoorig (x, u) = x, \quad \forall (x, u) \in \Xc \times \Uc,
\\
&\Eorigtoaug{u^*} (x) = (x, u^*), \quad \forall x \in \Xc , \quad
\forall u^* \in \Uc.  
\end{align}
\end{subequations}
Then,
\begin{enumerate}
\item  $\restr{f}{\Xc, u \equiv u^*} = f \circ \Eorigtoaug{u^*}$, for all
$f \in \Faug$ and $u^* \in \Uc$; 
\item  $g_e = g \circ \Raugtoorig$, for all $g \in \Fc$;
\item $\Raugtoorig \circ \Eorigtoaug{u^*} = \identity_{\Xc}$, for all $u^* \in \Uc$. 
\oprocend
\end{enumerate}
\end{lemma}

The proof
follows from the definitions.  We rely on
Lemma~\ref{l:mapping-based-augmented-original} to define appropriate
linear operators to move between the function spaces $\Fc$
and~$\Faug$.

\begin{proposition}\longthmtitle{Linear Operator Connection Between
$\Fc$ and $\Faug$}\label{p:linear-connection-augmented-original}
Assume $\Fc$ and $\Faug$ satisfy: 
\begin{enumerate}[label=($\mathfrak{C}$\roman*)]
\item $\restr{f}{\Xc, u \equiv u^*} \in \Fc$, for all $f \in \Faug$
and all $u^* \in \Uc$;
\item $g_e \in \Faug$, for all $g \in \Fc$.
\end{enumerate}
Define the family of operators
$\{\RopFaugtoF{u^*} : \Faug \to \Fc\}_{u^* \in \Uc}$ and the
operator $\EopFtoFaug : \Fc \to \Faug$ as
\begin{align*}
\RopFaugtoF{u^*} f = f \circ \Eorigtoaug{u^*} , \qquad \EopFtoFaug g = g \circ \Raugtoorig,
\end{align*}
for each $u^* \in \Uc$.
Then,
\begin{enumerate}
\item $\{\RopFaugtoF{u^*}\}_{u^* \in \Uc}$ and $\EopFtoFaug$ are well-defined and linear;
\item $\RopFaugtoF{u^*} f = \restr{f}{\Xc, u \equiv u^*}$, for all $f \in \Faug$ and $u^* \in \Uc$;
\item $\EopFtoFaug g = g_e$, for all $g \in \Fc$;
\item $\RopFaugtoF{u^*} \EopFtoFaug = \identity_{\Fc}$, for all $u^* \in \Uc$;
\item $\EopFtoFaug(\Fc) = \Faugci$;
\item $\RopFaugtoF{u^*}(\Faug) = \RopFaugtoF{u^*}(\Faugci)= \Fc$, for
all $u^* \in \Uc$.
\end{enumerate}
\end{proposition}
\begin{proof}
(a) Given any $u^* \in \Uc$, for $f \in \Faug$, using Lemma~\ref{l:mapping-based-augmented-original}, we have
$\RopFaugtoF{u^*} f = f \circ \Eorigtoaug{u^*} = \restr{f}{\Xc, u
\equiv u^*}$. Therefore, $\RopFaugtoF{u^*} f \in \Fc$ and the
operator is well-defined. A similar argument holds for $\EopFtoFaug$. 
Linearity of the operators directly
follows from their definition.

(b)-(c) The proof follows from
Lemma~\ref{l:mapping-based-augmented-original}.

(d) Given $u^* \in \Uc$, for $g \in \Fc$, using Lemma~\ref{l:mapping-based-augmented-original}(c), one can write
$\RopFaugtoF{u^*} \EopFtoFaug g = g \circ \Raugtoorig \circ
\Eorigtoaug{u^*} = g \circ \identity_{\Xc} = g$.

(e) We first prove $\EopFtoFaug(\Fc) \subseteq \Faugci$. Let
$g \in \Fc$. By part~(c), $\EopFtoFaug g = g_e \in \Faug$. Moreover, by definition, $g_e$ is control-independent.
Hence, $\EopFtoFaug g \in \Faugci$.

Next, we prove $\Faugci \subseteq \EopFtoFaug(\Fc)$. 
Given $f \in \Faugci$, select an arbitrary $u^* \in \Uc$ and consider $g_f = \restr{f}{\Xc, u \equiv u^*} \in \Fc$. Note
\begin{align*}
[\EopFtoFaug (g_f)] (x,u) & = \restr{f}{\Xc, u \equiv u^*}(x) = f(x,u^*) 
\\
& = f(x,u), \quad \forall x \in \Xc, \, \forall u, u^* \in \Uc,
\end{align*}
where in the last equality we use that $f$ is control-independent. This proves $\Faugci \subseteq \EopFtoFaug(\Fc)$, and we conclude $\Faugci = \EopFtoFaug(\Fc)$.

(f) $ \RopFaugtoF{u^*}(\Faugci)= \Fc$ follows from (d) and (e). The
rest follows from the fact that
$\RopFaugtoF{u^*}(\Faug) \subseteq \Fc$ by definition of
$\RopFaugtoF{u^*}$ and the fact that
$\Fc =\RopFaugtoF{u^*}(\Faugci) \subseteq \RopFaugtoF{u^*}(\Faug)$ by
virtue of $\Faugci \subseteq \Faug$.
\end{proof}

\subsection{Implications for Control-Independent Functions}\label{s:implications-CI}

As the discussion so far has illustrated, control-independent functions
play a key role in the technical treatment. This is because the domain of functions in spaces $\Finf$ and $\Faug$ are different from the state space of the control system~\eqref{eq:control-system} and, therefore, one has to rely on control-independent functions to connect the action of the operators $\Kc^{\infty}$ and $\Kaug$ to the trajectories of~\eqref{eq:control-system}, cf. Lemmas~\ref{l:capture-system-control-ind-inf} and~\ref{l:capture-system-control-ind-aug}, resp. These observations warrant a closer study of control-independent functions.

We first show that under reasonable conditions, there is a bijective relationship between $\Faugci$ and $\Finfci$.

\begin{proposition}\longthmtitle{Isomorphism\footnote{An isomorphism between two \emph{abstract} vector spaces is a bijective \emph{linear} map. Note that the linearity implies that it preserves the structure of the vector space and the two vector spaces are essentially the same under linear operations.} Between $\Faugci$ and $\Finfci$}\label{p:operator-connection-ci-inf-aug}
Assume $\Faug$ and $\Fc^\infty$ satisfy ($\mathtt{C}$i)-($\mathtt{C}$ii) in Proposition~\ref{p:linear-connection-infinite-augmented}. Define the extension $\EopCIFaugtoFinf: \Faugci \to \Finfci$ and restriction $\RopCIFinftoFaug:  \Finfci \to \Faugci$ operators as
\begin{align*}
\EopCIFaugtoFinf f = \EopFaugtoFinf f, \qquad \RopCIFinftoFaug g = \RopFinftoFaug g .
\end{align*}
Then,
\begin{enumerate}
\item $\EopCIFaugtoFinf$ and $\RopCIFinftoFaug$ are well-defined;
\item $\EopCIFaugtoFinf  \RopCIFinftoFaug = \identity_{\Finfci}$;
\item $\RopCIFinftoFaug  \EopCIFaugtoFinf = \identity_{\Faugci}$;
\item $\EopCIFaugtoFinf$ and $\RopCIFinftoFaug$ are bijective.
\end{enumerate}
\end{proposition}
\begin{proof}
(a) This follows from combining Proposition~\ref{p:linear-connection-infinite-augmented}(a), (e), and (f).

(b) Note that since operators have matching domain and codomain, $\EopCIFaugtoFinf  \RopCIFinftoFaug$ is well-defined. For $g \in \Finfci$, we have
\begin{align*}
\EopCIFaugtoFinf  \RopCIFinftoFaug g(x,\mathbf{u})  = g \circ \Eaugtoinf \circ \Rinftoaug (x, \mathbf{u})
\\
= g \circ \Eaugtoinf (x, \mathbf{u}(0)) = g(x, (\mathbf{u}(0), \mathbf{u}(0),\ldots)).
\end{align*}
Since $g$ is control-independent, $g(x, (\mathbf{u}(0), \mathbf{u}(0),\ldots)) = g(x,\mathbf{u})$ for all $(x, \mathbf{u})  \in \Xc \times \ell(\Uc)$, and therefore $\EopCIFaugtoFinf  \RopCIFinftoFaug g = g$.

(c) This follows from Proposition~\ref{p:linear-connection-infinite-augmented}(d).

(d) This is a consequence of~(b) and (c).
\end{proof}

Proposition~\ref{p:operator-connection-ci-inf-aug} has a major difference with respect to Proposition~\ref{p:linear-connection-infinite-augmented}: when the domain and codomain of operators are restricted to control-independent functions, the restriction and extensions operators between $\Faugci$ and $\Finfci$ become bijective and are inverse of each other. This means that, for every control-independent function in $f \in \Faugci$, there is a \emph{unique} control-independent function in $g \in \Finfci$ with $ g = \EopCIFaugtoFinf f$ such that their values on all states match,
\begin{align}\label{eq:matching}
f(x,u) = g(x,\mathbf{u}), \quad \forall x \in \Xc, \; \forall u \in \Uc, \; \forall \mathbf{u} \in \ell(\Uc).
\end{align}
The reverse of this statement also holds: for all $g \in \Finfci$, there is a unique $f \in \Faugci$ with $f = \RopCIFinftoFaug g$ such that~\eqref{eq:matching} holds.

We are also interested in connecting the spaces $\Faugci$ and~$\Fc$. We first show that the family of restrictions from $\Faug$ to $\Fc$ coincide when their domain is restricted to~$\Faugci$.

\begin{lemma}\longthmtitle{Restriction Operators from $\Faug$ to $\Fc$ Coincide on $\Faugci$}\label{l:restriction-coincide-Faugci}
Assume $\Fc$ and $\Faug$ satisfy ($\mathfrak{C}$i)-($\mathfrak{C}$ii) in Proposition~\ref{p:linear-connection-augmented-original}.
Then
\begin{align*}
\RopFaugtoF{u_1} f = \RopFaugtoF{u_2} f, \quad \forall f \in \Faugci, \; \forall u_1, u_2 \in \Uc.
\end{align*}
\end{lemma}
\begin{proof}
This follows by noting that, for $f \in \Faugci$, we have $
[\RopFaugtoF{u_1} f](x) = f(x,u_1) = f(x,u_2) =  [\RopFaugtoF{u_2} f](x)$, for all $x \in \Xc$, and all $u_1, u_2 \in \Uc$.
\end{proof}

As a result of Lemma~\ref{l:restriction-coincide-Faugci}, we define the \emph{restriction operator}
$\RopCIFaugtoF: \Faugci \to \Fc$ as
\begin{align*}
\RopCIFaugtoF f = \RopFaugtoF{u^*} f, 
\end{align*}
where $u^* \in \Uc$ is arbitrary. Next, we fully connect $\Faugci$ and $\Fc$ via bijective maps.

\begin{proposition}\longthmtitle{Isomorphism Between $\Faugci$ and $\Fc$}\label{p:operator-connection-ci-aug-orig}
Assume $\Fc$ and $\Faug$ satisfy ($\mathfrak{C}$i)-($\mathfrak{C}$ii) in Proposition~\ref{p:linear-connection-augmented-original}.
Define the \emph{extension operator} $\EopCIFtoFaug: \Fc \to \Faugci$  as
\begin{align*}
\EopCIFtoFaug f = \EopFtoFaug f .
\end{align*}
Then,
\begin{enumerate}
\item $\EopCIFtoFaug$ is well-defined;
\item $\EopCIFtoFaug  \RopCIFaugtoF = \identity_{\Faugci}$;
\item $\RopCIFaugtoF  \EopCIFtoFaug = \identity_{\Fc}$;
\item $\EopCIFtoFaug$ and $\RopCIFaugtoF$ are bijective.
\end{enumerate}
\end{proposition}
\begin{proof}
(a) This follows from Proposition~\ref{p:linear-connection-augmented-original}(a) and~(e). 

(b) For $g \in \Faugci$, we have
\begin{align*}
\EopCIFtoFaug  \RopCIFaugtoF g (x,u) = g \circ \Eorigtoaug{u^*} \circ \Raugtoorig (x,u) 
\\
= g(x,u^*) = g(x,u), \forall x \in \Xc, \; \forall u,u^*\in \Uc,
\end{align*}
where in the last equality we have used the fact that $g$ is control-independent.

(c) This follows from Proposition~\ref{p:linear-connection-augmented-original}(d).

(d) This is a consequence of (c) and (d).
\end{proof}

Proposition~\ref{p:operator-connection-ci-aug-orig} shows that there is a one-to-one correspondence between $\Fc$ and $\Faugci$. One  can think of $\Faugci$ as a copy of $\Fc$, with the domain changed from $\Xc$ to $\Xc \times \Uc$. For example, if $f \in \Fc$, then there exist $g = \EopCIFtoFaug f \in \Faugci$ where 
\begin{align}\label{eq:matching2}
f(x) = g(x,u), \; \forall x \in \Xc, \; \forall u \in \Uc.
\end{align}
The inverse of this statement is also true: for all $g \in \Faugci$ there is a unique $f= \RopCIFaugtoF g\in \Fc$ for which~\eqref{eq:matching2} holds.

The next result is a consequence of Propositions~\ref{p:operator-connection-ci-inf-aug} and~\ref{p:operator-connection-ci-aug-orig}.

\begin{corollary}\longthmtitle{Isomorphism Between $\Fc$, $\Faugci$, $\Finfci$}\label{c:isomorphism-control-invariant}
Assume $\Fc$, $\Faug$, and $\Fc^\infty$ satisfy
($\mathtt{C}$i)-($\mathtt{C}$ii) in
Proposition~\ref{p:linear-connection-infinite-augmented} and
($\mathfrak{C}$i)-($\mathfrak{C}$ii) in
Proposition~\ref{p:linear-connection-augmented-original}. Then,
\begin{enumerate}
\item $\Fc$, $\Faugci$, and $\Finfci$ are isomorphic;
\item $\card(\Fc) = \card(\Faugci) = \card(\Finfci)$;
\item $\card(\Fc) \leq  \card( \Faug ) \leq \card( \Finf)$.
\end{enumerate}
\end{corollary}
\begin{proof}
(a)-(b) follow from Propositions~\ref{p:operator-connection-ci-inf-aug} and~\ref{p:operator-connection-ci-aug-orig}. Regarding (c), the inequality $\card(\Fc) \leq \card(\Faug)$ follows from~(b) and the fact that $ \Faugci  \subseteq \Faug$. The other inequality follows from the fact that the operator $\EopFaugtoFinf: \Faug \to \Finf$ is injective since it has a left inverse based, cf. Proposition~\ref{p:linear-connection-infinite-augmented}(d).
\end{proof}

\begin{figure}[htb]
\centering 
{\includegraphics[width=.99\linewidth]{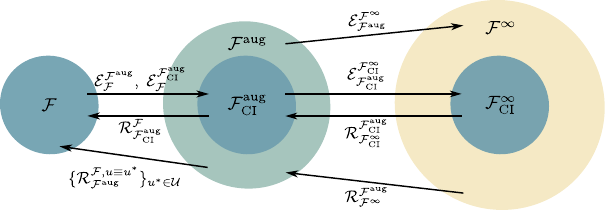}}
\caption{Connections between the functions spaces. The action of $\EopFtoFaug$ and $\EopCIFtoFaug$ coincide even though the have different codomains.}\label{fig:maps-function-spaces}
\vspace*{-1ex}
\end{figure}

Corollary~\ref{c:isomorphism-control-invariant} ensures that the spaces $\Fc$, $\Faug$, $\Finf$ have the same level of richness in terms of encoding the information about the trajectories of control system~\eqref{eq:control-system}. However, one should keep in mind that the spaces $\Faug$ and $\Finf$ generally  have larger cardinality than $\Fc$, cf. Figure~\ref{fig:maps-function-spaces}, since the frameworks $\Kaug$ and $\Kc^{\infty}$  embed  the effect of the input sequence in the function spaces while in the KCF framework the input information is embedded in the switching signal (and kept separate from~$\Fc$).

\section{Connections Between the Associated Koopman Operators}\label{sec:equivalence}
In this section, we use the augmented Koopman operator as an intermediary to
connect the Koopman operators associated with the infinite input
sequences framework and the Koopman Control Family. We perform this
process at two levels: dynamics and operators acting on
function spaces, cf. Table~\ref{tab:framework-comparison}.

\subsection{Infinite Input Sequences Framework and Augmented Koopman Operator}

We are ready to relate $\Tc^\infty$ and $\Taug$ using the notions introduced in Section~\ref{s:functions-inf-aug}.

\begin{proposition}\longthmtitle{Dynamics Connection Between
Infinite Input Sequences and Augmented Koopman Operator
Frameworks}\label{p:system-connection-infinite-augmented}
Let $\Eaugtoinf$ and $\Rinftoaug$ be the maps defined
in~\eqref{eq:mapping-based-infinite-augmented}. Then, 
\begin{enumerate}
\item $\Tc^\infty \circ \Eaugtoinf = \Eaugtoinf \circ \Taug$;
\item $\Taug = \Rinftoaug \circ \Tc^\infty \circ \Eaugtoinf$. 
\end{enumerate}
\end{proposition}
\begin{proof}
(a) For $(x, u) \in \Xc \times \Uc$, we have
\begin{align*}
\Tc^\infty \circ \Eaugtoinf (x, u) & = \Tc^\infty (x, (u, u, \ldots))
\nonumber\\
& = (\Tc(x,u), (u,u, \ldots)).
\end{align*}
On the other hand, 
\begin{align*}
\Eaugtoinf \circ \Taug (x, u) = \Eaugtoinf (\Tc(x,u), u)
\nonumber\\
= (\Tc(x,u), (u,u, \ldots)).
\end{align*}
Therefore, $\Tc^\infty \circ \Eaugtoinf (x, u) = \Eaugtoinf \circ \Taug (x, u)$.

(b) For $(x, u) \in \Xc \times \Uc$, we have
\begin{align*}
&\Rinftoaug \! \circ \! \Tc^\infty \! \circ \! \Eaugtoinf (x, u) 
\!=\! \Rinftoaug \!\circ \! \Tc^\infty (x, (u, u, \ldots))
\nonumber\\
&= \Rinftoaug (\Tc(x,u), (u,u, \ldots)) 
\!=\! (\Tc(x,u), u) \!=\! \Taug (x, u) .
\end{align*}
\end{proof}

Proposition~\ref{p:system-connection-infinite-augmented} provides a
tool based on domain restriction and extension maps to connect the
action of dynamical systems in the infinite input sequences and the
augmented Koopman operator frameworks. The
connections are via function compositions, allowing to state
equivalent descriptions at the operator level.

\begin{theorem}\longthmtitle{Operator Connection Between
Infinite Input Sequences and Augmented Koopman Operator
Frameworks}\label{t:operator-connection-infinite-augmented}
Assume $\Faug$ and $\Fc^\infty$ satisfy ($\mathtt{C}$i)-($\mathtt{C}$ii) in Proposition~\ref{p:linear-connection-infinite-augmented}. 
Then,
\begin{enumerate}
\item $\RopFinftoFaug \Kc^\infty = \Kaug \RopFinftoFaug$;
\item $\Kaug = \RopFinftoFaug \Kc^\infty \EopFaugtoFinf$;
\item $\restr{(\Kc^\infty f)}{\Xc \times \Uc} = \Kaug (\restr{ f}{\Xc \times \Uc})$, for all $f \in \Fc^\infty$;
\item $\Kaug g = \restr{(\Kc^\infty g^\infty)}{\Xc \times \Uc}$, for all $g \in \Faug$. 
\end{enumerate}
\end{theorem}
\begin{proof}
(a) For $h \in \Fc^\infty$, we have
\begin{align*}
\RopFinftoFaug \Kc^\infty h & = \RopFinftoFaug (h \circ \Tc^\infty) = h \circ \Tc^\infty \circ \Eaugtoinf.
\end{align*}
This, together with 
Proposition~\ref{p:system-connection-infinite-augmented}(a), yields
\begin{align*}
h \circ \Tc^{\infty} \circ \Eaugtoinf
&= h \circ \Eaugtoinf \circ \Taug 
\\
&= \Kaug (h \circ \Eaugtoinf)
\\
&= \Kaug (\RopFinftoFaug h ).
\end{align*}
Therefore, $\RopFinftoFaug \Kc^\infty h = \Kaug \RopFinftoFaug h$.

(b) For $g \in \Faug$, we have
\begin{align*}
\RopFinftoFaug \Kc^\infty \EopFaugtoFinf g
&= \RopFinftoFaug \Kc^\infty (g \circ \Rinftoaug)
\nonumber \\
&=\RopFinftoFaug( g \circ \Rinftoaug \circ \Tc^{\infty}) 
\nonumber \\
&= g \circ \Rinftoaug \circ \Tc^\infty \circ \Eaugtoinf.
\end{align*}
This, together with
Proposition~\ref{p:system-connection-infinite-augmented}(b), yields
\begin{align*}
\RopFinftoFaug \Kc^\infty \EopFaugtoFinf g = g \circ \Taug = \Kaug g.
\end{align*}

(c)-(d) This follows from~(a)-(b) in conjunction with
Proposition~\ref{p:linear-connection-infinite-augmented}.
\end{proof}

\begin{figure}[htbp]
\centering
\adjustbox{scale=1.2,center}{%
\begin{tikzcd}[scale=10.4, row sep=large, column sep=large]
\Faug \arrow[r, "\EopFaugtoFinf"] \arrow[d, "\Kaug"'] \arrow[rr, "\identity_{\Faug}", bend left=40] 
& \Finf \arrow[r, "\RopFinftoFaug"] \arrow[d, "\Kc^{\infty}"] & \Faug \arrow[d, "\Kaug"]  \\
\Faug                   & \Finf \arrow[l, "\RopFinftoFaug"]   \arrow[r, "\RopFinftoFaug" ']                & \Faug 
\end{tikzcd}
}
\caption{Commutative diagram illustrating
Theorem~\ref{t:operator-connection-infinite-augmented}.}
\label{fig:commutative-aug-inf}
\end{figure}
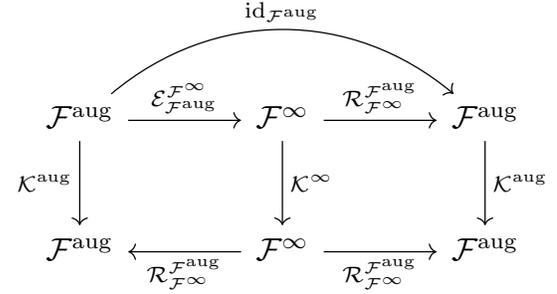

Figure~\ref{fig:commutative-aug-inf} shows the commutative diagram for
the operators' actions described in
Theorem~\ref{t:operator-connection-infinite-augmented}.

\subsection{Koopman Control Family and Augmented Koopman Operator}

We are ready to relate $\{\Tc_{u^*}\}_{u^* \in \Uc}$ and $\Taug$ using the notions introduced in Section~\ref{s:functions-ori-aug}.

\begin{proposition}\longthmtitle{Dynamics Connection Between
Koopman Control Family and Augmented Koopman Operator
Frameworks}\label{p:system-connection-augmented-KCF}
Let $\Raugtoorig$ and $\{\Eorigtoaug{u^*}\}_{u^* \in \Uc}$ be the
maps defined in~\eqref{eq:maps-e-r-aug-orig}. Then, for all $u^* \in \Uc$,
\begin{enumerate}
\item
$\Taug \circ \Eorigtoaug{u^*} = \Eorigtoaug{u^*} \circ \Tc_{u^*}$;
\item $\Tc_{u^*} = \Raugtoorig \circ \Taug \circ \Eorigtoaug{u^*} $. 
\end{enumerate}
\end{proposition}
\begin{proof}
Let $u^* \in \Uc$. (a) for $x \in \Xc$, we have
\begin{align*}
\Taug \circ \Eorigtoaug{u^*} (x) = \Taug (x, u^*) = (\Tc(x,u^*), u^*).
\end{align*}
On the other hand,
\begin{align*}
\Eorigtoaug{u^*} \circ \Tc_{u^*} (x) = \Eorigtoaug{u^*}
(\Tc(x,u^*)) = (\Tc(x,u^*), u^*). 
\end{align*}
Therefore, $\Taug \circ \Eorigtoaug{u^*} = \Eorigtoaug{u^*} \circ \Tc_{u^*}$.

(b) for $x \in \Xc$, we have
\begin{align*}
\Raugtoorig  \circ \Taug \circ \Eorigtoaug{u^*} (x) = \Raugtoorig \circ \Taug (x, u^*)  
\nonumber
\\
= \Raugtoorig  (\Tc(x,u^*), u^*) = \Tc(x,u^*) = \Tc_{u^*} (x).
\end{align*}
Therefore, $\Raugtoorig  \circ \Taug \circ \Eorigtoaug{u^*} = \Tc_{u^*}$.
\end{proof}

This result allows us to state equivalent descriptions at the operator
level between the augmented Koopman operator and the Koopman control family frameworks. 

\begin{theorem}\longthmtitle{Operator Connection Between Koopman
Control Family and Augmented Koopman Operator
Frameworks}\label{t:operator-connection-augmented-KCF}
Assume $\Fc$ and $\Faug$ satisfy ($\mathfrak{C}$i)-($\mathfrak{C}$ii) in Proposition~\ref{p:linear-connection-augmented-original}.
Then, for all
$u^* \in \Uc$
\begin{enumerate}
\item $\RopFaugtoF{u^*} \Kaug = \Kc_{u^*} \RopFaugtoF{u^*}$;
\item $\Kc_{u^*} = \RopFaugtoF{u^*} \Kaug \EopFtoFaug$;
\item $\restr{(\Kaug f)}{\Xc, u \equiv u^*} = \Kc_{u^*} (\restr{ f}{\Xc,
u \equiv u^*})$, for all $f \in \Faug$;
\item $\Kc_{u^*} g = \restr{(\Kaug g_e)}{\Xc, u \equiv u^*}$, for all $g \in \Fc$. 
\end{enumerate}
\end{theorem}
\begin{proof}
Let
$u^* \in \Uc$. (a) For $f \in \Faug$, we have
\begin{align*}
\RopFaugtoF{u^*} \Kaug f \!=\! \RopFaugtoF{u^*} (f \circ \Taug) 
\!=\! f \circ \Taug \circ \Eorigtoaug{u^*}.
\end{align*}
On the other hand, 
\begin{align*}
\Kc_{u^*} \RopFaugtoF{u^*} f \!=\! \Kc_{u^*} (f \circ \Eorigtoaug{u^*})
\!=\! f \circ \Eorigtoaug{u^*} \circ \Tc_{u^*}. 
\end{align*}
Using Proposition~\ref{p:system-connection-augmented-KCF}(a), we can write 
$\Kc_{u^*} \RopFaugtoF{u^*} f = f \circ \Taug \circ \Eorigtoaug{u^*}$, and the result follows.

(b) For $g \in \Fc$, we have
\begin{align*}
\RopFaugtoF{u^*} \Kaug \EopFtoFaug g
&= \RopFaugtoF{u^*} \Kaug (g \circ \Raugtoorig) 
\nonumber
\\ 
&= \RopFaugtoF{u^*} (g \circ \Raugtoorig \circ \Taug) 
\nonumber \\
&= g \circ \Raugtoorig \circ \Taug \circ \Eorigtoaug{u^*}.
\end{align*}
Using
Proposition~\ref{p:system-connection-augmented-KCF}(b), we
deduce that
$\RopFaugtoF{u^*} \Kaug \EopFtoFaug g = g \circ \Tc_{u^*} = \Kc_{u^*}
g$.

(c)-(d)  This follows from (a)-(b) in conjunction with Lemma~\ref{l:mapping-based-augmented-original}.
\end{proof}

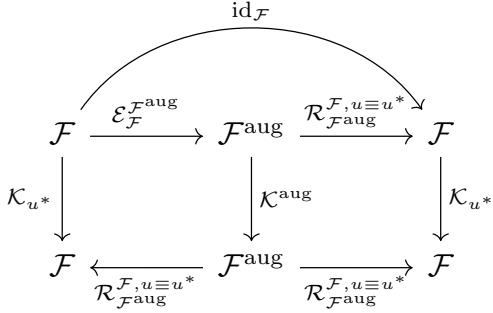
\begin{figure}[htbp]
\centering
\adjustbox{scale=1.2,center}{%
\begin{tikzcd}[scale=10.4, row sep=large, column sep=large]
\Fc 
\arrow[r, "\EopFtoFaug"] 
\arrow[d, "\Kc_{u^*}"'] 
\arrow[rr, "\identity_{\Fc}", bend left=55] 
& \Faug 
\arrow[r, "\RopFaugtoF{u^*}"] 
\arrow[d, "\Kaug"] 
& \Fc 
\arrow[d,"\Kc_{u^*}"]
\\
\Fc 
& \Faug 
\arrow[l, "\RopFaugtoF{u^*}"]   
\arrow[r, "\RopFaugtoF{u^*}" ']                
& \Fc 
\end{tikzcd}
}
\caption{Commutative diagram illustrating
Theorem~\ref{t:operator-connection-augmented-KCF}. Here, $u^* \in \Uc$ is arbitrary.}
\label{fig:commutative-aug-orig}
\end{figure}

Figure~\ref{fig:commutative-aug-orig} shows the commutative diagram
for the operators' actions described in
Theorem~\ref{t:operator-connection-augmented-KCF}.

\subsection{Equivalence Results on Control-Independent Functions}

Following up on Section~\ref{s:implications-CI}, here we turn our attention to the operators' action on control-independent functions. The next result reveals how the actions of $\{\Kc_{u^*}\}_{u^*\in \Uc}$, $\Kaug$, and $\Kinf$ are connected on the isomorphic spaces $\Fc$, $\Faugci$, and~$\Finfci$.

\begin{theorem}\longthmtitle{Operators' Actions on Control-Independent Functions}\label{t:operators-on-control-independent}
Assume $\Fc$, $\Faug$, and $\Fc^\infty$ satisfy
($\mathtt{C}$i)-($\mathtt{C}$ii) in
Proposition~\ref{p:linear-connection-infinite-augmented} and
($\mathfrak{C}$i)-($\mathfrak{C}$ii) in
Proposition~\ref{p:linear-connection-augmented-original}. Then,
\begin{enumerate}
\item $\Kc_{u^*} f = \RopFaugtoF{u^*} \Kaug \, \EopCIFtoFaug f$, for all $f \in \Fc$ and $u^* \in \Uc$;
\item $ \Kaug g = \RopFinftoFaug \Kc^{\infty} \EopCIFaugtoFinf g$,  for all $g \in \Faugci$;
\item $\Kc_{u^*} f =  \RopFaugtoF{u^*} \RopFinftoFaug \Kc^{\infty} \EopCIFaugtoFinf \EopCIFtoFaug f$, for all $ f \in \Fc$ and $u^* \in \Uc$.
\end{enumerate}
\begin{proof}
(a) This follows from Theorem~\ref{t:operator-connection-augmented-KCF}(b) and the definition of~$\EopCIFtoFaug$.

(b) This follows from Theorem~\ref{t:operator-connection-infinite-augmented}(b) and the definition of~$\EopCIFaugtoFinf$.

(c) This follows from (a) and (b).
\end{proof}
\end{theorem}

\begin{figure}[htbp]
\centering
\adjustbox{scale=1.2,center}{%
\begin{tikzcd}[row sep=large, column sep=large]
\Fc \arrow[d, "\Kc_{u^*}"] 
\arrow[r, "\EopCIFtoFaug", shift left] 
& \Faugci
\arrow[l, "\RopCIFaugtoF", shift left] 
\arrow[d, "\Kaug"] 
\arrow[r, "\EopCIFaugtoFinf", shift left] 
& \Finfci
\arrow[l, "\RopCIFinftoFaug", shift left] 
\arrow[d, "\Kc^{\infty}"] 
\\
\Fc 
& \Faug 
\arrow[l, "\RopFaugtoF{u^*}"] 
& \Finf
\arrow[l, "\RopFinftoFaug"] 
\end{tikzcd}
}
\caption{Commutative diagram illustrating
Theorem~\ref{t:operators-on-control-independent} and Propositions~\ref{p:operator-connection-ci-inf-aug} and~\ref{p:operator-connection-ci-aug-orig}. Here, $u^* \in \Uc$ is arbitrary.}
\label{fig:commutative-ci-orig-aug-inf}
\end{figure}
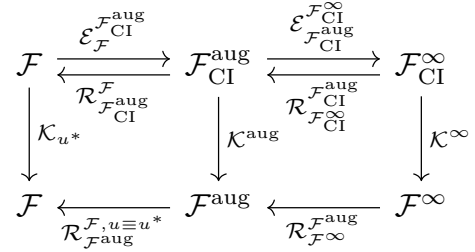

Figure~\ref{fig:commutative-ci-orig-aug-inf} illustrates the action of operators in Theorem~\ref{t:operators-on-control-independent} on different function spaces starting from control-independent functions. 

Figure~\ref{fig:commutative-ci-orig-aug-inf} reveals that one can start from any of the spaces $\Fc$, $\Faugci$, or $\Finfci$ through equivalent functions linked by isomorphisms $\EopCIFtoFaug$, $\RopCIFaugtoF$, $\EopCIFaugtoFinf$, and $\RopCIFinftoFaug$ and capture dynamical information from system~\eqref{eq:control-system} by applying $\{\Kc_{u^*}\}_{u^* \in \Uc}$, $\Kaug$, or $\Kinf$. Note that, after applying $\Kaug$ and $\Kinf$, the output functions will not be necessarily control-independent anymore, since the system behavior depends on the input. It also worth mentioning that the output functions in $\Faug$ and $\Finf$ capture both the information of the system's trajectories as well as the information of the input sequence (since the domain of functions in $\Faug$ and $\Finf$ are $\Xc \times \Uc$ and $\Xc \times \ell(\Uc)$, resp.). Therefore, the operators $\RopFinftoFaug$ and $\{\RopFaugtoF{u^*}\}_{u^* \in \Uc}$ act as filters, detaching the information about the input sequence except the first element influencing the systems' trajectory for the next time step. 

Next, we focus our attention on a different kind of equivalence result pertaining to the evolution of function values on multi-step trajectories under the Koopman operator via infinite input sequences and the Koopman Control Family framework. 

\begin{theorem}\longthmtitle{Equivalence of KCF and Koopman Operator via Infinite Input Sequences on Trajectories}\label{t:KCF-Kinf-trajectories}
Assume $\Fc$, $\Faug$, and $\Fc^\infty$ satisfy
($\mathtt{C}$i)-($\mathtt{C}$ii) in
Proposition~\ref{p:linear-connection-infinite-augmented} and
($\mathfrak{C}$i)-($\mathfrak{C}$ii) in
Proposition~\ref{p:linear-connection-augmented-original}. Let
$\{x_k\}_{k \in \naturals_0}$ be the trajectory
of~\eqref{eq:control-system} from initial condition $x_0$ with input
sequence $\mathbf{u} = (u_0, u_1, \ldots)$. Then, for all
$k \in \naturals_0$,
\begin{enumerate}
\item for all $f \in \Fc$, we have
\begin{align*}
f(x_k) &= [\Kc_{u_0} \Kc_{u_1} \ldots \Kc_{u_{k-1}} f] (x_0) 
\\
&= [(\Kinf)^k \EopCIFaugtoFinf \EopCIFtoFaug f] (x_0, \mathbf{u});
\end{align*}
\item for all $h \in \Finfci$ with decomposition $h = h_\Xc \, 1_{\ell(\Uc)}$ (cf. Definition~\ref{def:control-independent-func-inf}), we have
\begin{align*}
h_\Xc (x_k) &= [(\Kinf)^k h] (x_0,\mathbf{u})
\\
&= [\Kc_{u_0} \Kc_{u_1} \ldots \Kc_{u_{k-1}} \RopCIFaugtoF \RopCIFinftoFaug h](x_0).
\end{align*}
\end{enumerate}
\end{theorem}
\begin{proof}
(a) Given $f \in \Fc$,  the first equality follows from Lemma~\ref{l:encoding-trajectories-KCF}. Let $h := \EopCIFaugtoFinf \EopCIFtoFaug f \in \Finfci$, cf. Propositions~\ref{p:operator-connection-ci-inf-aug} and~\ref{p:operator-connection-ci-aug-orig}.  Since $h$ is control-independent, it can be decomposed as $ h = h_\Xc 1_{\ell(\Uc)}$. Based on Lemma~\ref{l:capture-system-control-ind-inf}, we have $h_\Xc (x_k) = [(\Kinf)^k h] (x_0,\mathbf{u})$. Hence, to prove the result, we seek to establish that $h_\Xc = f$. Based on the definition of the extension operators $\EopCIFtoFaug$ and $\EopCIFaugtoFinf$, one can write
\begin{align*}
h (x, \mathbf{\bar{u}}) &= [\EopCIFaugtoFinf \EopCIFtoFaug f] (x, \mathbf{\bar{u}}) 
\\
&= [f \circ \Raugtoorig \circ \Rinftoaug] (x, \mathbf{\bar{u}}) = [f \circ \Raugtoorig](x,\mathbf{\bar{u}}(0))
\\
& = f(x), \quad \forall x \in \Xc, \quad \forall \mathbf{\bar{u}} \in \ell(\Uc),
\end{align*}
and therefore $h_\Xc = f$. 

(b) Given $h \in \Finfci$, let $f:= \RopCIFaugtoF \RopCIFinftoFaug h \in \Fc$. Note that $h  = \EopCIFaugtoFinf \EopCIFtoFaug f$ based on Propositions~\ref{p:operator-connection-ci-inf-aug} and~\ref{p:operator-connection-ci-aug-orig}. The result now follows from (a).
\end{proof}

Theorem~\ref{t:KCF-Kinf-trajectories} provides a direct connection
between the KCF and the Koopman operator via infinite-sequences on the
system
trajectories.
It is important to point out that this result does not include the
augmented Koopman operator, $\Kaug$, since it is not a Koopman
operator associated with the control system~\eqref{eq:control-system}
and does not directly capture multi-step trajectories.

We finish this section by explaining a fundamental difference between
Theorems~\ref{t:operators-on-control-independent}
and~\ref{t:KCF-Kinf-trajectories}. Theorem~\ref{t:operators-on-control-independent}
pertains to the operators' action on function spaces, revealing how
the functions encode the system's information in general settings, and
how one can separate the information of system trajectories from the
information of input sequences and change the domain of functions to
the state space of control system~\eqref{eq:control-system}. On the
other hand, Theorem~\ref{t:KCF-Kinf-trajectories} reveals information
about the function \emph{values} on trajectories as opposed to the
functions themselves. This result is more useful in modeling and
control applications while
Theorem~\ref{t:operators-on-control-independent} pertains to deeper
fundamental connections between the different frameworks.

\section{Conclusions}
We have studied the connections between two extensions of Koopman
operator theory to control systems, the Koopman operator via infinite
input sequences and the Koopman control family. Since each extension
relies on a different mechanism to encode system information into
operators acting on vector spaces, we first examined how the
information of trajectories and input sequences is captured in
each. This understanding then enabled us to provide ways to connect
the function spaces via linear composition operators.  We relied on
these operators as bridges to connect the actions of Koopman-based
formulations in each framework. As a result, we provided a
comprehensive analysis of their structure, along with constructive
algebraic recipes to convert the actions of the operators. Finally, we
showed that, under mild conditions on the function spaces, the
frameworks are equivalent, both in terms of encoding system
information in function spaces and the evolution of function values
along system trajectories. Future work will explore the implications
of these connections for studying nonlinear control systems as well as
for control design procedures.

\bibliographystyle{IEEEtran}%
\bibliography{short}

\vspace*{-5ex}

\begin{IEEEbiography}[{\includegraphics[width=1in,height=1.2in,clip,keepaspectratio]{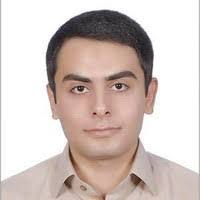}}]{Masih Haseli} 
received the B.Sc. and M.Sc. degrees in Electrical Engineering from Amirkabir University of Technology (Tehran Polytechnic), Tehran, Iran, in 2013 and 2015, respectively. He obtained the Ph.D. degree in Engineering Sciences (Mechanical Engineering) from the University of California, San Diego, CA, USA, in 2022. From 2022 to 2025, he held a postdoctoral position in the Department of Mechanical and Aerospace Engineering at the University of California, San Diego. He is currently a postdoctoral scholar at the California Institute of Technology, Pasadena, CA, USA. His research interests include system identification, nonlinear systems, network systems, data-driven modeling and control, and distributed and parallel computing. Dr. Haseli is the recipient of the Bronze Medal in the 2014 Iran National Mathematics Competition and the Best Student Paper Award at the 2021 American Control Conference.
\end{IEEEbiography}

\vspace*{-3ex}

\begin{IEEEbiography}[{\includegraphics[width=1.1in,keepaspectratio]{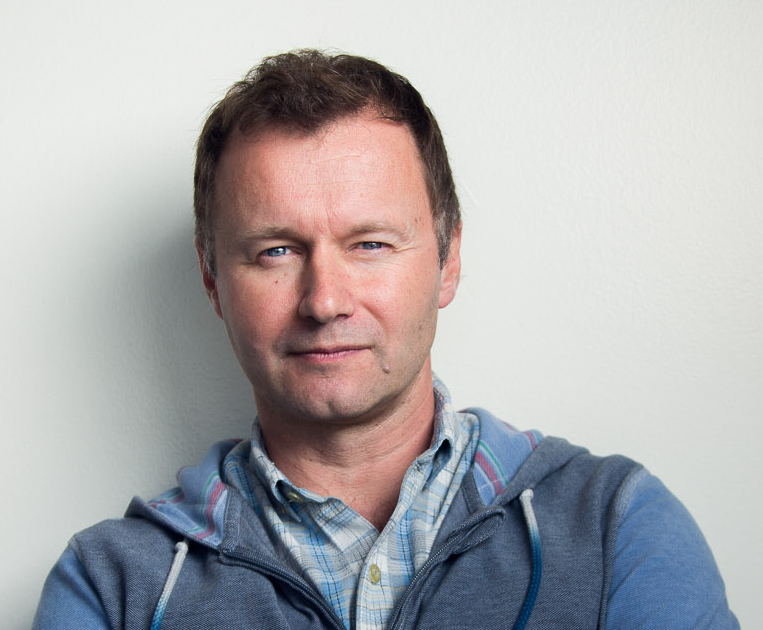}}]{Igor Mezi\'c} 
(M'08,SM'20,F'22) received the Dipl. Ing. in mechanical
engineering in 1990 from the University
of Rijeka, Rijeka, Croatia, and the Ph.D. degree in
applied mechanics from the California Institute of
Technology.
He was a postdoctoral researcher at the Mathematics
Institute, University of Warwick, Warwick,
U.K., during 1994–1995. From 1995 to 1999, he was
a member of the Mechanical Engineering Department,
University of California Santa Barbara, where
he is currently a Distinguished Professor. In 2000–2001, he was an
Associate Professor in the Division of Engineering and Applied Sciences at
Harvard University, Cambridge, MA. His work  in the fields of dynamical
systems
and control theory centers around the applications and advances in
theory of Koopman operators.
Dr. Mezic was awarded the Alfred P. Sloan Fellowship, the NSF CAREER
Award from the NSF, the George S. Axelby Outstanding Paper Award on
“Control of Mixing” from the IEEE, and the United Technologies Senior Vice
President for Science and Technology Special Achievement Prize in 2007.
He is a fellow of the IEEE, SIAM and APS. He received the biennial SIAM
J.D Crawford Prize in Dynamical Systems in 2021.
He was an Editor of Physica D: Nonlinear Phenomena and an Associate
Editor of the Journal of Applied Mechanics and SIAM Journal on Control
and Optimization. He is the Director of the Center for Energy Efficient
Design
and Head of Buildings and Design Solutions Group at the Institute for Energy
Efficiency at the University of California at Santa Barbara.
\end{IEEEbiography}

\vspace*{-3ex}

\begin{IEEEbiography}[{\includegraphics[width=1in,height=1.2in,clip,keepaspectratio]{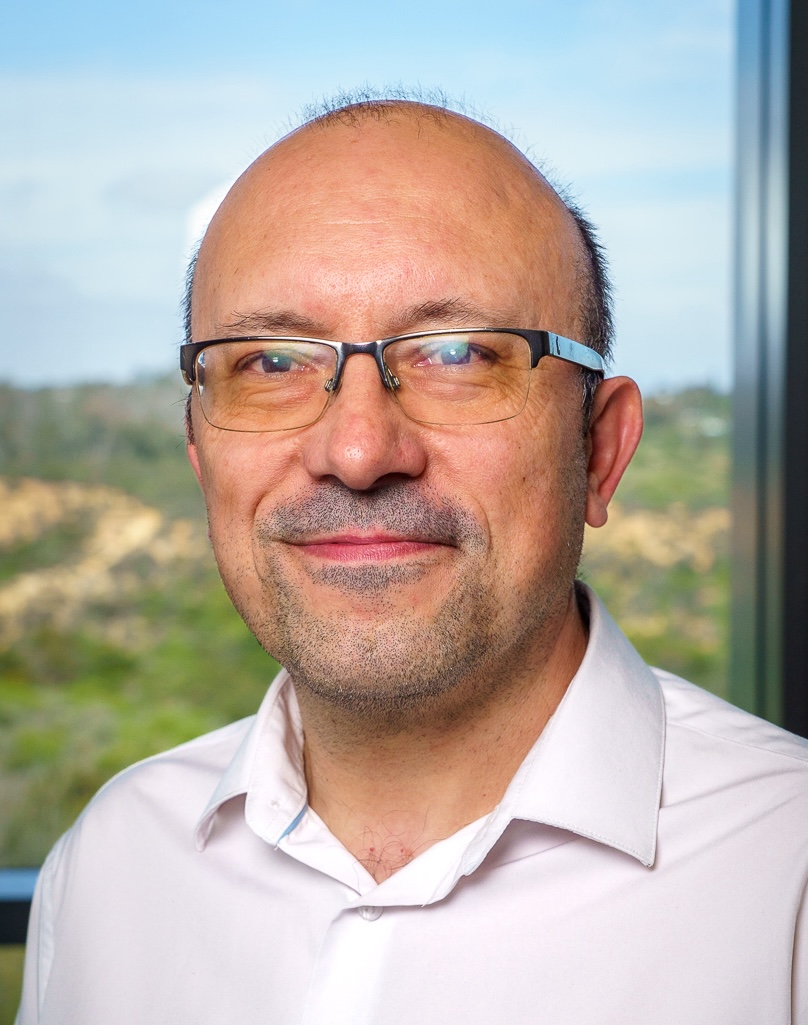}}]{Jorge
Cort\'{e}s}(M'02, SM'06, F'14) received the Licenciatura degree in
mathematics from Universidad de Zaragoza, Spain, in 1997, and the
Ph.D. degree in engineering mathematics from Universidad Carlos III
de Madrid, Spain, in 2001. He held postdoctoral positions with the
University of Twente, Twente, The Netherlands, and the University of
Illinois at Urbana-Champaign, Illinois, USA.
He is a Professor and Cymer Corporation Endowed Chair in High
Performance Dynamic Systems Modeling and Control at the Department
of Mechanical and Aerospace Engineering, UC San Diego, California,
USA.  He is a Fellow of IEEE, SIAM, and IFAC.  His research
interests include distributed control and optimization, network
science, nonsmooth analysis, reasoning and decision making under
uncertainty, network neuroscience, and multi-agent coordination in
robotic, power, and transportation networks.
\end{IEEEbiography}

\end{document}